%% file: paper.tex
\documentclass[a4paper]{article}

\usepackage[a4paper, total={5.125in, 8.5in}]{geometry}
\usepackage{hyperref}
\usepackage[utf8]{inputenc}
\usepackage[english]{babel}
\usepackage{cleveref}
\usepackage{amsmath,amsfonts}
\usepackage{mathabx}
\usepackage{textgreek}
\usepackage{bbm}
\usepackage[vlined,boxruled,linesnumbered,algo2e,algosection]{algorithm2e}
\usepackage{xcolor}
\usepackage{tikz}
\usetikzlibrary{intersections,plotmarks}

\usepackage{pgfplotstable,filecontents}
\pgfplotsset{compat=1.14}
\usepackage{csquotes}
\usepackage{enumerate}
\usepackage{graphicx}
\graphicspath{{./figs/}}

\newcommand{\SetAlgorithmStyle}{
  \setcounter{AlgoLine}{0}
  \SetKwData{Left}{left}\SetKwData{This}{this}\SetKwData{Up}{up}
  \SetKwInOut{Input}{Input}
  \SetKwInOut{Output}{Output}
  \ResetInOut{input}
  \SetKwComment{tcp}{//}{}
  \SetKwFor{For}{for}{}{end}
  \SetArgSty{}
  \DontPrintSemicolon
}

\input{macros.tex}

\crefalias{AlgoLine}{line}
\crefalias{algocf}{algorithm}

\makeatletter
\let\cref@old@stepcounter\stepcounter
\def\stepcounter#1{%
  \cref@old@stepcounter{#1}%
  \cref@constructprefix{#1}{\cref@result}%
  \@ifundefined{cref@#1@alias}%
    {\def\@tempa{#1}}%
    {\def\@tempa{\csname cref@#1@alias\endcsname}}%
  \protected@edef\cref@currentlabel{%
    [\@tempa][\arabic{#1}][\cref@result]%
    \csname p@#1\endcsname\csname the#1\endcsname}}
\makeatother

\makeatletter
\renewcommand{\algocf@caption@boxruled}{%
  \hrule
  \hbox to \hsize{%
    \vrule\hskip-0.4pt
    \vbox{   
       \vskip\interspacetitleboxruled%
       \unhbox\algocf@capbox\hfill
       \vskip\interspacetitleboxruled
       }%
     \hskip-0.4pt\vrule%
   }\nointerlineskip%
}%
\makeatother

\title{Coarsening in Algebraic Multigrid using Gaussian Processes\thanks{This work was partially funded by Deutsche Forschungsgemeinschaft (DFG) Transregional Collaborative Research Centre 55 (SFB/TRR55)}}
\author{Hanno Gottschalk\thanks{Bergische Universiat\"at Wuppertal, IMACM, Faculty of Mathematics and Natural Sciences, \texttt{\{hgotsch,kkahl\}@uni-wuppertal.de}}\ \and Karsten Kahl\footnotemark[2]}

\begin{document}
    \maketitle
    
    \begin{abstract}
    Multigrid methods have proven to be an invaluable tool to efficiently solve large sparse linear systems arising in the discretization of partial differential equations (PDEs). Algebraic multigrid methods and in particular adaptive algebraic multigrid approaches have shown that multigrid efficiency can be obtained without having to resort to properties of the PDE. Yet the required setup of these methods poses a not negligible overhead cost. 
    Methods from machine learning have attracted attention to streamline processes based on statistical models being trained on the available data.
    Interpreting algebraically smooth error
    as an instance of a Gaussian process, we develop a new, data driven approach to construct adaptive algebraic multigrid methods. Based on Gaussian a priori distributions,  Kriging interpolation minimizes the mean squared error of the a posteriori distribution, given the data on the coarse grid. Going one step further, we exploit the quantification of uncertainty in the Gaussian process model in order to construct efficient variable splittings. Using a semivariogram fit of a suitable covariance model we demonstrate that our approach yields efficient methods using a single algebraically smooth vector.
    \end{abstract}
    
    \section{Introduction}\label{sec:intro}
    The solution of large sparse linear systems of equations
    \begin{equation}\label{eq:linsys}
    Ax = b,\ b\in \mathbb{R}^{n},\ A \in \mathbb{R}^{n \times n} \text{\ symmetric positive definite,}
    \end{equation} that arise in the discretization of partial differential equations, typically makes up the bulk of computations in modern scientific computing. It is thus of utmost importance to come up with efficient algorithms to solve these systems of equations. By exploiting a separation of scales multigrid methods can achieve optimal linear complexity for this task, but heavily rely of the availability of expert knowledge about the particular partial differential equation that the linear systems originates from as well as the employed discretization scheme; cf.~\cite{BrenScot2008}.
    Due to the fact that this information might not readily be available or that there is no known geometric multigrid construction, the concept of algebraic multigrid methods has been introduced in~\cite{Bran1986,BranMcCoRuge1985,RugeStue1987,Stue1983,XuZika2017}. Efficiency in algebraic multigrid methods is achieved by pairing a simple iterative scheme, the \textit{smoother}, with a variational coarse grid correction. Assuming a smoother is defined by $M \approx A^{-1}$, the error propagator of a two-grid algebraic multigrid method with Galerkin coarse grid construction is given by
    \begin{equation} \label{eq:eprop}
        E_{2g} = (I-MA) (I-P (P^T A P)^{-1} P^T A) (I-MA).
    \end{equation} Due to the variational construction of coarse grid correction the whole setup of an algebraic multigrid method can be reduced to the definition of the interpolation operator $P$. In particular, the dimension of the coarse space, $n_{c}$, the interpolation relations, i.e., the sparsity pattern of $P\in \mathbb{R}^{n\times n_{c}}$ and its entries need to be defined. Typically these tasks are split into two parts. Finding $n_{c}$ and the sparsity pattern of $P$ is often referred to as the \textit{coarsening} problem, while determining the entries of $P$ is known as the \textit{interpolation} problem.
    
    In the first algebraic multigrid methods~\cite{BranMcCoRuge1985,RugeStue1987} operator based approaches have been suggested to solve both problems. In case $A$ has $M$-matrix structure, e.g., as a particular discretization of an elliptic partial differential equation, it can be shown that these approaches lead to methods with fast convergence. However, these early approaches rely heavily on assumptions about the underlying problem and therefore cannot be extended significantly beyond the $M$-matrix case. A huge step in overcoming this limitation has been the introduction of adaptivity in algebraic multigrid methods~\cite{BranBranKahlLivs2011,BrezFalgMacLMantMcCoRuge2006,BrezFalgMacLMantMcCoRuge2004}. 
    Common in all adaptive approaches in algebraic multigrid methods is the idea to guide the definition of interpolation, posed in terms of the coarsening and interpolation problem either by using spectral information about $A$ and/or the smoothing iteration. Due to the fact that explicit calculation of (partial) spectra is prohibitively expensive these methods rely on an iterative approximation process that makes use of the emerging multigrid hierarchy. In contrast to the interpolation problem, where several approaches showed promosing performance~\cite{BranBranKahlLivs2015b,BranBranKahlLivs2011,MandBrezVane1999,OlsoSchrTumi2011,WanChanSmit1999,XuZika2004}, the coarsening problem turned out to be harder to tackle.
    
    Many of the approaches that have been tried to solve the coarsening problem in adaptive algebraic multigrid methods revolve around the idea of strength of connection, a concept introduced in classical algebraic multigrid. This includes approaches based on binary variable relations such as~\cite{BranBranKahlLivs2015b,BranChenKrauZika2013,BranChenZika2012,OlsoSchrTumi2011} but approaches that take relations of more than two variables into account such as~\cite{KahlRott2018}.
    The detection of strongly connected pairs of variables is also found in aggregation-based approaches such as~\cite{BranChenKrauZika2013,LivnBran2012,NapoNota2016,Nota2010}.
    In some sense the idea of compatible relaxation, introduced in~\cite{BranFalg2010}, comes closest to general applicability, but is hard to integrate and mesh with typical solutions to the interpolation problems, i.e., the definition of the entries of $P$, in adaptive algebraic multigrid approaches.
    
    In this paper we propose a new way of solving the coarsening problem which resembles the least squares interpolation approach of the bootstrap algebraic multigrid framework~\cite{BranBranKahlLivs2011,KahlRott2018}. By considering the test vectors of the bootstrap framework as instances of a Gaussian process~\cite{adler2010geometry,bogachev1998gaussian,vanmarcke2010random} we are able to apply techniques from machine learning, 
    especially the concept of conditional (a posteriori) distributions. To calibrate this statistical model to the data, we use parametric
     semivariogram  models to fit the covariance structure to the data provided by algebraically smooth test vectors. 
     This enables us to efficiently  solve the coarsening and interpolation problem at the same time. The idea of Kriging interpolation is to view the values of test vectors at coarse grid variables as partial observations of a Gaussian process. Based on these observed values at the coarse grid variables, expected values and errors at the fine grid variables can be derived by  computing conditional expectations and variances. This method  originally stems from spatial interpolation in geostatistics, see e.g. \cite{christakos2012random} 
     but has been widely used in various machine learning and engineering tasks in the past \cite{forrester2008engineering,kleijnen2009kriging}. Its main advantage lies in its statistical properties as best linear unbiased predictor, given the data of the field on the 'observed locations' -- the coarse grid -- and the spatial correlation structure of the data. Coincidentially, the Kriging iterpolator bears some resemblence to the least-squares formulation of interpolation introduced in~\cite{BranBranKahlLivs2011} and more specifically the operator-based modification of it found in~\cite{MantMcCoParkRuge2010}.
    
    Inheriting from the spatial correlation structure of the Gaussian process, the conditional variance given the observations on certain points, is low at points with sufficient observations in the neighborhood and high elsewhere. This additional piece of information can now be used if one has to decide on where to make the next observation. In this paper we apply a greedy optimization procedure, picking the point of highest conditional variance given coarse grid results, to adaptively refine the coarse grid by subsequently adding fine grid points to the coarse grid until the conditional variance on all remaining fine grid points is small.

    A statistical view on the adaptive setup in algebraic multigrid methods is not completely new. It has been used in~\cite{LivnBran2012} to motivate a definition of strength of connection based on a measure of correlation present in test vectors, but this construction lacks the framing of Gaussian processes.
    Other related work in the context of the integration of ideas from stochastics into multigrid method has been  presented in~\cite{owhadi2017multigrid}. While this work also contains a deep mathematical analysis of convergence, there are several points where our work takes a different route. First of all, our work is purely algebraic and only uses structures that can be derived from the matrix $A$. Hence we do not use any structures that stem from the spatial structures of the underlying PDE and thus we do not have techniques at our disposal that use spectral equivalence and other techniques based on harmonic analysis. Instead, we only use distance notions between 'nodes' that can intrinsically be derived from $A$. Also, we resort to data driven estimation of correlation structures instead of an analytic derivation of these structures from the matrix $A$. This empirical approach largely reduces the computational cost in the choice of priors. 
    Furthermore, while the prior chosen in \cite{owhadi2017multigrid} is supported by error analysis in the $A$-norm, it has the disadvantage of producing a generalized random field with distributional paths \cite{gelfand1964generalized}, which is a disputable choice for the prior belief on the solution to a PDE. The price we pay is to use a more experimental and less deeply founded approach.     
    
    In~\cref{sec:multigrid} we give a short introduction into the construction of algebraic multigrid methods and highlight how adaptivity can be used in order to capture the nature and underlying structure of the problem at hand, especially in terms finding suitable coarsenings. Then we give an overview on Gaussian processes in~\cref{sec:gaussian}, where we motivate the connection to the adaptive setup process and discuss interpolation in the context of these processes. This leads us directly to the formulation of Kriging interpolation in~\cref{sec:linear}, which we discuss and compare to the least squares interpolation approach in~\cref{app:kriging-ls}. Finally, using additional heuristics we present our new adaptive coarsening approach in~\cref{sec:coarsening} before closing with numerical tests in~\cref{sec:numerics} demonstrating the potential of our approach and some final remarks in~\cref{sec:conclusion}.
    
    \section{Adaptive algebraic multigrid methods}\label{sec:multigrid}
\begin{figure}
        \begin{center}
            \begin{tikzpicture}
            \node at (0,0) {\includegraphics[width=.9\textwidth]{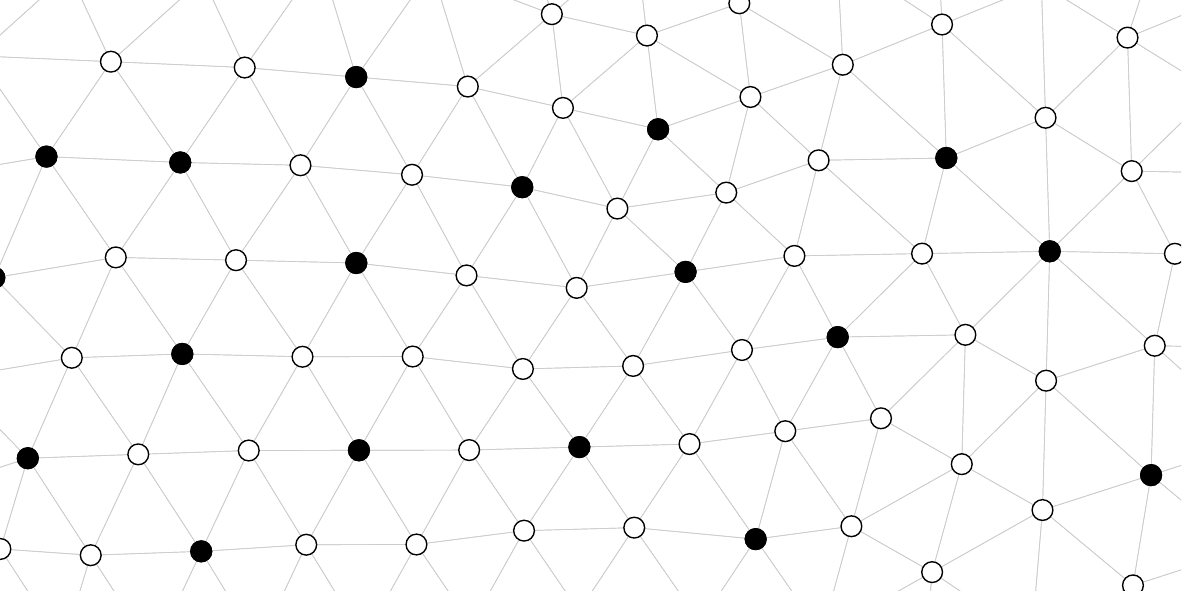}};
            \end{tikzpicture}
        \end{center}
        \caption[Splitting of the set of variables $\mathcal{V}$ into a set of coarse variables $\mathcal{C}$ and a set of fine variables  $\mathcal{F}$.]{Splitting of the set of all variables $\mathcal{V}$ into a set of coarse variables $\mathcal{C}$ depicted by \raisebox{-.1em}{\resizebox{.75em}{!}{\tikz{\node[circle,fill=black,minimum width=1em] at (0,0) { };}}}
        and a set of fine variables  $\mathcal{F}$ depicted by \raisebox{-.1em}{\resizebox{.75em}{!}{\tikz{\node[circle,fill=white,draw=black,minimum width=1em] at (0,0) { };}}}
        .}\label{fig:variable_splitting}
    \end{figure}

    The efficiency of multigrid methods lies in the complementarity of the smoothing iteration and the coarse grid correction. Algebraic multigrid methods construct complementarity without relying on knowledge of the underlying problem or the employed discretization strategy. Assuming that $A$ is symmetric positive definite, it is common to consider a Galerkin construction for the coarse grid correction error propagator
    \[
        E_{\rm cgc} = I - P \left(P^{T}AP\right)^{-1}P^{T}A.
    \] Thus it is completely determined by the definition of the interpolation operator $P: \mathbb{R}^{n_c} \rightarrow \mathbb{R}^{n}$. Due to the fact that in this case $E_{\rm cgc}$ corresponds to the $A$-orthogonal projection onto the space $A$-orthogonal to $\operatorname{range}(P)$, we can assume for simplicity sake and a better intuition of the construction of $P$, that it can be represented in the following form
    \begin{equation}\label{eq:canonical_interpolation}
        P = \begin{bmatrix} P_{\mathcal{F}\!,\mathcal{C}}\\ I \end{bmatrix}.
    \end{equation} This particular form of $P$ can be assumed due to the fact that $E_{\rm cgc}$ remains unchanged when transforming $P \rightarrow PX$ for any non-singular $X$, i.e., it depends solely on $\operatorname{range}(P)$ and not on the basis representation of this subspace. Using this form of $P$, we observe, that the index set of all variables $\mathcal{V} = \{1,2,\ldots,n\}$ can be split into two disjoint sets. The set $\mathcal{C}$ of variables that define the coarse grid and the remaining set $\mathcal{F} = \mathcal{V}\setminus\mathcal{C}$ of variables solely being present in the fine grid system as depicted in~\cref{fig:variable_splitting}. 
    
    In that sense~\eqref{eq:canonical_interpolation} defines interpolation from the variables with $\mathcal{C}$ indices, which are kept identical when moving from coarse to fine variables, to $\mathcal{F}$ variables with interpolation weights found in $P_{\rm \mathcal{F}\!,\mathcal{C}}$. As already mentioned in the introduction the definition of interpolation now reduces to three questions. First, which coarse variable set $\mathcal{C} \subset \mathcal{V}$ to choose, in particular this also amounts to determining the coarsening ratio $\tfrac{|\mathcal{C}|}{n}$. Second, for each variable $i$ in $\mathcal{F}$, a set of interpolatory variables $\mathcal{C}_i \subset \mathcal{C}$ has to be defined, which corresponds to the sparsity pattern of $P_{\rm \mathcal{F}\!,\mathcal{C}}$. Last, the entries of $P_{\rm \mathcal{F}\!,\mathcal{C}}$ have to be defined such that complementarity of the smoothing iteration and the coarse grid correction is achieved. In addition, as an implicit requirement, the sparsity of the coarse system of equations, given by $P^{T}AP$ has to be guaranteed in order to be able to apply the construction recursively and thus achieve optimal linear complexity.
    Based on the findings in~\cite{BranCaoKahlFalgHu2018} the complementarity of the smoothing iteration and the coarse-grid correction is equivalent to the requirement that $\operatorname{range}(P)$ approximates the space spanned by eigenvectors of the error propagator of the smoother corresponding to eigenvalues close to $1$, i.e., components that are slow to converge--also known as algebraically smooth error components~\cite{RugeStue1987}.
    
    The adaptive construction of algebraic multigrid methods can thus be interpreted as the generation of a low dimensional, sparse representation of the space of algebraically smooth error. In the context of the bootstrap algebraic multigrid framework and its extensions~\cite{BranBranKahlLivs2011,KahlRott2018,MantMcCoParkRuge2010} this is facilitated by the use of a set of algebraically smooth test vectors $\{v^{(1)},\ldots,v^{(K)}\}$ which are obtained by smoothing initially random test vectors with entries that stem from a normal distribution. By drawing on a connection to Gaussian processes we are now going to modify the construction of interpolation in this setting using ideas that originate in statistical geophysics.

    \section{Gaussian Processes}\label{sec:gaussian}
    In this section, we introduce Gaussian processes \cite{adler2010geometry} and the Kriging predictor \cite{bivand2008applied,sherman2011spatial} to interpolate from coarse to fine grid points. We discuss various approaches for modeling the covariance of the underlying Gaussian process, including semivariogram estimation based on test vectors.    
    
    \subsection{\label{sec:stochProc} Gaussian Stochastic Processes and Kriging}
    A stochastic process with the index set $\mathcal{I}$ is a collection of real valued random variables $X=\{X_i\mid i\in \mathcal{I}\}$ on a common probability space $(\Omega,\Sigma,\mathbb{P})$, where $\mathbb{E}[\cdot]$ denotes the expected value with respect to $\mathbb{P}$. Here $\Omega$ is the event set, $\Sigma$ the sigma field and $\mathbb{P}$ the probability measure. $X$ is a Gaussian process, if the distribution of any finite collection $X_\mathcal{H}=(X_{i_1},\ldots,X_{i_q})^T$ for finite subsets $\mathcal{H}=\{i_1,\ldots,i_q\}\subseteq \mathcal{I}$ is a multivariate Gaussian distribution $N(\mu_\mathcal{H},{C}_\mathcal{H})$ with probability density
    \begin{equation}
        \label{eq:densGauss}
        f_\mathcal{H}(x_\mathcal{H})=\frac{1}{\sqrt{2\pi}^q|{C}_\mathcal{H}|^{\frac{1}{2}}}e^{-\frac{1}{2}(x_\mathcal{H}-\mu_\mathcal{H})^T{C}_\mathcal{H}^{-1}(x_\mathcal{H}-\mu_\mathcal{H})},~~x_\mathcal{H}\in\mathbb{R}^m,
    \end{equation}
    where $|\cdot|$ denotes the determinant, $\mu_\mathcal{H}\in\mathbb{R}^q$ is the expected value and  ${C}_\mathcal{H}\in \mathbb{R}^{q\times q}$ is a positive definite covariance matrix. Consistency of the finite dimensional distributions in the sense of Kolmogrov \cite{tao2011introduction} implies that there exist functions 
    \[
    \mu:\mathcal{I}\to\mathbb{R} \text{\ and\ } {C}:\mathcal{I}\times\mathcal{I}\to \mathbb{R}
    \] such that $\mu_\mathcal{H}=(\mu_{i_1},\ldots,\mu_{i_q})^T$ and ${C}_\mathcal{H}=({C}_{i,j})_{i,j\in \mathcal{H}}$. 
    
    In the given context, we interpret $\mathcal{I}$ as the computational domain, and consider $\mathcal{V} \subseteq \mathcal{I}$ as a finite discretization of $\mathcal{I}$. $X_i$ represents the epistemic uncertainty about the solution of $\eqref{eq:linsys}$ at the grid point $i\in\mathcal{V}$. 
    
    Let us suppose we gathered partial information on $X_{\mathcal{V}}$, e.g., by solving $\eqref{eq:linsys}$ on a coarse subset of variables $\mathcal{C} \subseteq \mathcal{V}$ and we would like to infer about the solution on the whole grid $\mathcal{V}$. By construction, any prediction $\widehat X_i, i\in \mathcal{V}$ can only depend on information from $X_\mathcal{C}$ and therefore has to be measurable with respect to the sigma field $\sigma_\mathcal{C}\subseteq\Sigma$ 
    associated with $X_{\mathcal{C}}$. We thus wish to find an optimal solution to the problem of making a prediction based on the values of $X_{\mathcal{C}}$ that minimizes the expected squared error, also called mean square error (\emph{MSE}), for $X_i, i \in \mathcal{V}$
    \begin{equation}
        \label{eq:minError}
        \mathbb{E}\left[\left(X_i-\widehat X_i\right)^2\right]\to \min \text{\quad s.t.\quad } \widehat X_i \text{\ is\ } \sigma_\mathcal{C}\text{-measurable}.
    \end{equation}
    Let $L^2(\Omega,\Sigma,\mathbb{P})$ be the space of square integrable random variables, then $L^2(\Omega,\sigma_{\mathcal{C}},\mathbb{P})\linebreak\subseteq L^2(\Omega,\sigma,\mathbb{P})$ of $\sigma_\mathcal{C}$-measurable $L^2$ functions is a closed subspace. Thus, the problem \eqref{eq:minError} is uniquely solved by the conditional expected value $\widehat X_i \in \mathbb{E}[X_i|X_\mathcal{C}]$ defined as the $L^2$-projection of $X_i\in L^2(\Omega,\Sigma,\mathbb{P})$ to  $L^2(\Omega,\sigma_{\mathcal{C}},\mathbb{P})$. This definition immediately implies the interpolation property $\widehat X_i=X_i$ for all $i\in \mathcal{C}$ as $X_i$ in this case is $\sigma_{\mathcal C}$ measurable itself. 
    
    To calculate the conditional expected value $\mathbb{E}[X_i|X_\mathcal{C}]$ along with the minimum expected squared error \eqref{eq:minError} on the remaining set of variables $\mathcal{F} = \mathcal{V}\setminus\mathcal{C}$, we first compute the conditional distribution of $X_{\mathcal{F}}$ given $X_\mathcal{C} = x_\mathcal{C}$ via its density as
    \begin{align*}
    \begin{split}
    f_{\mathcal{F}|\mathcal{C}}(x_{\mathcal{F}}|x_\mathcal{C})
    & = \frac{f_{\mathcal{V}}\left(x_{\mathcal{F}},x_\mathcal{C}\right)}{f_\mathcal{C}(x_\mathcal{C})}\\
    & = \frac{1}{\sqrt{2\pi}^q|{C}_{\mathcal{F}|\mathcal{C}}|^{\frac{1}{2}}}e^{-\frac{1}{2}(x_{\mathcal{F}}-\mu_{\mathcal{F}}(x_\mathcal{C}))^T{C}_{\mathcal{F}|\mathcal{C}}^{-1}(x_{\mathcal{F}}-\mu_{\mathcal{F}}(x_{\mathcal{C}}))}, 
    \end{split}
    \end{align*}
    where from \eqref{eq:densGauss} we get from straight forward calculation
    \begin{subequations}
    \begin{align}\label{eq:kriging}
        \mu_{\mathcal{F}}(x_{\mathcal{C}}) 
        &= \mu_{\mathcal{F}}+{C}_{\mathcal{F}\!,\mathcal{C}}{C}_\mathcal{C}^{-1}\!\left(x_\mathcal{C}-\mu_{\mathcal{C}}\right),\\
        \label{eq:krigVar}
        {C}_{\mathcal{F}|\mathcal{C}} 
        &= {C}_{\mathcal{C}} - {C}_{\mathcal{F}\!,\mathcal{C}}{C}_{\mathcal{C}}^{-1}{C}_{\mathcal{F}\!,\mathcal{C}}^{T},
    \end{align}
    \end{subequations}
    where we introduced the notation ${C}_{\mathcal{H}\!,\mathcal{J}}=(C_{i,j})_{i\in \mathcal{H}, j \in \mathcal{J}}$ for any finite $\mathcal{H}\!,\mathcal{J}\subseteq \mathcal{I}$.
    
    Using \eqref{eq:kriging} and \eqref{eq:krigVar} we see that for any $i \in \mathcal{F}$ the prediction $\widehat X_i = \mu_{i}(X_{\mathcal{C}})$, i.e., the conditional expected value, minimizes the MSE and its conditional variance is given by
    \begin{align*}
    \begin{split}
        \sigma^2_{i|\mathcal{C}} &= {C}_{i}-{C}_{i,\mathcal{C}}{C}_{\mathcal{C}}^{-1}{C}_{i,\mathcal{C}}^{T}\\
        &= \min_{\widehat X_i
        }\mathbb{E}\left[\left(X_i-\widehat{X}_i\right)^2\right].
        \end{split}
    \end{align*}
    
    
    
    \subsection{\label{sec:linear}Linear Interpolation from Kriging}
    Note that \eqref{eq:kriging} provides an affine-linear interpolation rule and not a linear one as required for the construction of the matrix $P$ in \eqref{eq:canonical_interpolation}. This problem can be dealt with in two ways: First, we can set the expected value $\mu_\mathcal{\mathcal{V}}=0$. This is consistent with the estimation of the data $\mu_\mathcal{V}$ and $C_\mathcal{V}$ defining the Gaussian process on the entire fine grid $\mathcal{V}$ from test vectors, see subsection \ref{sec:covMod}. Alternatively, assuming that $\sigma_\mathcal{V}$ is fixed or already estimated and that $\mu_{i}\cong \mu$ is constant, we obtain the value of $\mu$ as the best linear unbiased predictor \emph{BLUP} from the data $x_\mathcal{C}$ on $\mathcal{C}$. In fact, suppose that we estimate $\mu$ linearly by $\widehat \mu=w_{\mathcal{C}}^T x_{\mathcal{C}}$. The requirement that this estimator is unbiased results in
    \begin{align*}
    \begin{split}
        \mathbb{E}_\mu\left[w^T_\mathcal{C}X_\mathcal{C}\right] &= w^T_\mathcal{C}\mathbb{E}_\mu\left[X_\mathcal{C}\right]=w^T_\mathcal{C}\mathbbm{1}_\mathcal{C}\mu\\
        &\Rightarrow\ w^T_\mathcal{C}\mathbbm{1}_\mathcal{C}=1,
    \end{split}
    \end{align*}
    where $\mathbbm{1}_\mathcal{C}$ is the vector of all ones, i.e., $(\mathbbm{1}_{\mathcal{C}})_i=1, i\in\mathcal{C}$ and $\mathbb{E}_\mu$ stands for the expected value for the Gaussian process with constant mean $\mu$. The optimal set of weights $w_{\mathcal{C}}$, given $C_\mathcal{V}$ and hence by restriction $C_\mathcal{C}$, is obtained by the following constrained optimization problem
    \begin{equation}\label{eq:krigingMeana}
    w_\mathcal{C} \in \operatorname{argmin} \left\{\mathbb{E}_\mu\left[(w^TX_\mathcal{C}-\mu)^2\right],\ w^T\mathbbm{1}_\mathcal{C}=1\right\}.
    \end{equation} Using a standard Lagrangian approach, we can reformulate this constrained quadratic optimization problem to the following set of equations
    \begin{equation}\label{eq:krigingMeanb}
        \frac{\partial L(w,\lambda)}{\partial w}=0 \text{\ and\ } \frac{\partial L(w,\lambda)}{\partial \lambda}=0,
    \end{equation} where $L(w,\lambda) = w^TC_\mathcal{C}w-\lambda(w^T\mathbbm{1}_\mathcal{C}-1)$. 
    It is now easily seen that the solution $w_{\mathcal{C}}$ of the equations \eqref{eq:krigingMeanb} is given by
    \begin{equation*}
        w = \frac{\mathcal{C}^{-1}_\mathcal{C}\mathbbm{1}_\mathcal{C}}{\mathbbm{1}_\mathcal{C}^TC_\mathcal{C}^{-1}\mathbbm{1}_\mathcal{C}^T}\ \Longrightarrow\ \widehat\mu=\frac{\mathbbm{1}_\mathcal{C}^T\mathcal{C}^{-1}_\mathcal{C}X_\mathcal{C}}{\mathbbm{1}_\mathcal{C}^TC_\mathcal{C}^{-1}\mathbbm{1}_\mathcal{C}^T}
    \end{equation*}
    Inserting $\mu_\mathcal{C}=\frac{\mathbbm{1}_\mathcal{C}^T\mathcal{C}^{-1}_\mathcal{C}x_\mathcal{C}}{\mathbbm{1}_\mathcal{C}^TC_\mathcal{C}^{-1}\mathbbm{1}_\mathcal{C}^T}\mathbbm{1}_{\mathcal{C}}$ into \eqref{eq:kriging} results in a prediction that depends linearly on the observed data.
    Thus this predictor corresponds to the construction of an interpolation matrix $P=P_{\mathcal{F}\!,\mathcal{C}}$ in \eqref{eq:canonical_interpolation}. 
    
    From both alternatives described here, we follow the second, estimating the value of $\mu$ from the coarse grid data $x_\mathcal{C}$ and not setting it to zero. The reason is that, while for test vectors $\mu=0$ is a natural choice (see subsection \ref{sec:covMod}), this is not necessarily the case for the problem, the multigrid solver is applied to.  A more thorough comparison of both approaches goes beyond the scope of this initial study. Note, that so far, we have only considered the calculation of the entries of $P_{\mathcal{F}\!,\mathcal{C}}$ given a subset $\mathcal{C}$ and neglecting the sparsity requirement of $P$ for now. The choice of $\mathcal{C}$ and the construction of localized, i.e., sparse, interpolation will be discussed next.
  
    \subsection{\label{sec:locKrig} Local Kriging}
    The computational cost of \eqref{eq:kriging} and \eqref{eq:krigVar} in many cases is prohibitive, due to the fact that $C_\mathcal{C}$ in general is not sparse. However, we can localize~\eqref{eq:kriging} and~\eqref{eq:krigVar} in the following sense. Assuming that there exists a (pseudo) distance $d_\mathcal{V}(i,j)$ on $\mathcal{V}$ and that the correlation 
    \[
        \varrho(i,j)=\frac{C_{i,j}}{\sqrt{{C}_{i,i}{C}_{j,j}}}
    \] decreases sufficiently fast when $d_\mathcal{V}(i,j)$ grows, we can neglect the effect of observations in far away points in $\mathcal{C}$ on the prediction $\widehat X_i$ of $X_i$.  In many cases, it is therefore sufficient to choose a subset $\mathcal{C}_i$ of $\mathcal{C}$ containing the $q_{\rm max}$ points in $\mathcal{C}$ that are closest to $i\in \mathcal{F}$ for the calculation of a suitable predictor. The number of neighbours $q_{\rm max}$ can be chosen independently of the size of $\mathcal{C}$ and $\mathcal{V}$ and is referred to as the \emph{caliber} of interpolation.
    
    We note that \eqref{eq:kriging} and \eqref{eq:krigVar} remain valid when replacing $\mathcal{C}$ with $\mathcal{C}_i$ and $\mathcal{F}$ with $i$. Also, the estimates \eqref{eq:krigingMeana} and \eqref{eq:krigingMeanb} can be localized accordingly. The complexity of computing the predictions on $\mathcal{F}$ therefore acquires the optimal linear growth in the size of this set, provided the construction $\mathcal{C}_i$ is either negligible for practical purposes or can be carried out with optimal complexity as well. Due to the applied localization this can be guaranteed by employing efficient graph based techniques. 
    
    
    \subsection{\label{sec:covMod} Non-Parametric Covariance Estimation}
    The selection of the expected value $\mu_\mathcal{V}$ and the covariance $C_\mathcal{V}$ has to essentially capture the correlation structure of the problem at hand. Assuming that a set of test vectors $V = \begin{bmatrix}v^{(1)}&\mid &\cdots&\mid&v^{(K)}\end{bmatrix}$ is given on $\mathcal{V}$, we define the average value and the empirical covariance matrix of these test vectors by  
    \begin{subequations}
    \begin{align}
    \label{eq:average}
        \widehat{\mu}_\mathcal{V}&=\frac{1}{K}V\mathbbm{1}_K,\\
        \label{eq:empCov}
        \widehat C&= \frac{1}{K}\left(V-\frac{1}{K}V\mathbbm{1}_K\mathbbm{1}_K^T\right) \left(V-\frac{1}{K}V\mathbbm{1}_K\mathbbm{1}_K^T\right)^T ,
    \end{align}
    \end{subequations}
    where $\mathbbm{1}_K\in\mathbb{R}^K$ is the column vector with value $1$ in every entry. By construction the rank of the empirical covariance matrix is bounded from above by $K$, the number of test vectors, so that it is in general not possible to simply replace the theoretical covariance $C$ in \eqref{eq:kriging} and \eqref{eq:krigVar} by the empirical counterparts. As long as $K<|\mathcal{C}|$, i.e., the number of test vectors is smaller than the number of coarse grid variables, $\widehat{C}_\mathcal{C}^{-1}$ does not exist. 
    There are two ways to fix this problem. First, we can regularize the calculation by replacing $\widehat C_\mathcal{C}\to \widehat C_\mathcal{C}+\varepsilon I$ where $\varepsilon > 0$ stands for a substitution $X\to X+N$ where $N=\{N_i\}, i\in \mathcal{I}$ is spatially uncorrelated Gaussian white noise with zero mean and variance $\varepsilon$.
    Second, in the context of local Kriging described in subsection \ref{sec:locKrig} it is only required that $\widehat C_{\mathcal{C}_i}$ is positive definite for all $i\in\mathcal{F}$. As we can choose $|C_{\mathcal{C}_i}| = q_{\rm max}\leq K$, this guarantees the non-singularity of $\widehat C_{\mathcal{C}_i}$ for almost all sets of randomly generated testvectors $V$ so that in this case an $\varepsilon$-regularisation is not required.

    Note that in many cases, the test vectors $v^{(j)}$ will be constructed by the application of one or a few smoothing, e.g., Gauss-Seidel, steps applied to a vector with noise data on $\mathcal{V}$, which is statistically centred around $0$. In such cases, $\widehat\mu_{\mathcal{V}}=0$ by theoretical considerations and we can replace the statistical estimation in \eqref{eq:average}. Likewise, in this case we could as well simplify \eqref{eq:empCov} by omission of the terms $\frac{1}{K}V\mathbbm{1}_K\mathbbm{1}_K^T$. 
    
     There exist several choices for the distance function in local Kriging. The first option is to use the coordinate distance $d^c(i,j)$, if an embedding of $\mathcal{V}$ in $\mathbb{R}^d$ is known and the underlying continuum problem is isotropic. This will not always be the case, especially in the context of algebraic multigrid. Instead we use the graph distance $d^A(i,j)$, which measures the shortest path in the undirected graph associated with the system matrix $A$ of~\eqref{eq:linsys} over $\mathcal{V}$ assuming that the length of an edge is defined as the inverse of the corresponding matrix entry, i.e., edge $\{i,j\}$ has length $\frac{1}{|A_{i,j}|}$. 
    
    In the context of local Kriging, it is not necessary to assemble the full matrix $\widehat C$, but only $n_{\mathcal{F}} = |\mathcal{F}|$ submatrices  $\widehat C_{\mathcal{C}_{i}}$ of size $q_{\rm max}\times q_{\rm max}$  and $n_{\mathcal{F}}$ matrices $\widehat C_{i,\mathcal{C}_{i}}$ of size $1\times q_{\rm max}$. Each entry of these matrices requires flops  proportional to $K$, which gives linear complexity in $n_{\mathcal{F}}$, provided the search for $\mathcal{C}_{i}$ is either negligible in terms of compute time on the relevant problem sizes or is implemented with optimal complexity in the sense that the search for $\mathcal{C}_{i}$ has bounded complexity. 
    \subsection{Parametric Semivariogram Estimation}\label{sec:semiVar}
    The disadvantage in the procedure described above lies in the fact that a relatively large number $K$ of test vectors $v^{(j)}$ is required in order to obtain a reasonable estimate $\widehat C$ for the underlying covariance structure $C$. Given that each test vector is calculated with a computational cost proportional to $n$, the generation of up to $K\approx 100$ test vectors can be a significant computational burden.
    
    Also, the large number $n$ of estimates of matrices $\widehat C_{\mathcal{C}_i}$ increases the probability that there is at least one $i\in\mathcal{\mathcal{F}}$ for which the estimate for $\widehat C_{\mathcal{C}_i}$ is poor, making it more difficult to obtain guarantees for the estimation of $\widehat X_i$ for all $i\in\mathcal{F}$.
    
    To reduce $K$ down to numbers in the range $1$ -- $10$ and to stabilise the estimation of single elements in $C_{\mathcal{C}_i}$, we follow methods from geostatistics that allow, under appropriate assumptions, a more efficient estimate.  
    
    To this purpose, we assume that the underlying continuous problem has some kind of translation invariance.  This can either be caused by a strict invariance of the underlying operator which is discretized by $A$ (neglecting the effect of boundary conditions) or an invariance in some statistical sense, where the local inhomogeneity is statistically the same around all points in the computational domain. In both cases it is legitimate, to work with translation invariant models for the a priori distribution of the Gaussian process.
    
    A stochastic process $X=\{X_i\}, i\in \mathcal{I}$ with the special choice $\mathcal{I}=\mathbb{R}^d$ is stationary, if the process $X_h=\{X_{i+h} \mid i\in \mathcal{I}\}$ has the same finite dimensional distributions as $X$ for all $h\in\mathbb{R}^d$. By \eqref{eq:densGauss} for Gaussian processes this amounts to a constant $\mu_i$ independent of $i$ and $C_{i+h,j+h}=C_{i,j}$ for all $i,j,h\in\mathbb{R}^d$. Consequently, the covariance function $C_{i,j} = C_{0,i-j} = :C(i-j)$ defines a function $C:\mathbb{R}^d\to \mathbb{R}$. Furthermore, the process is called isotropic, if $C(h)=C(\Lambda h)$ for any rotation matrix $\Lambda\in SO(d)$ and $h\in \mathbb{R}^d$. This can be justified if the underlying operator in the continuum (approximately) shares this property. In this case, we can model (with a slight abuse of notation) $C(h)=C(|h|)$ where only a function $C:\mathbb{R}_+\to\mathbb{R}$ has to be estimated from the test vectors $V$. This is a standard task in geostatistics \cite{bivand2008applied,sherman2011spatial}, which we briefly review.
    
    In geostatistics it is customary to fit semivariograms instead of the covariance function. Both are connected via
    \begin{align}
    \begin{split}\label{eq:semiVar}
        \gamma(|h|)&=\frac{1}{2}\mathbb{E}\left[(X(0)-X(h))^2\right]\\
        &=\frac{1}{2}\mathbb{E}\left[(X(0)-\mu)^2+2(X(0)-\mu)(X(h)-\mu)+(X(h)-\mu)^2\right]\\
        &=C(0)-C(|h|)
        \end{split}
    \end{align}
    and contain the same information, as the asymptotic sill value $C(0)$ can be obtained by letting $h\to\infty$ and thereby $C(|h|)\to 0$ as correlations decline at large distances. Note, that in \eqref{eq:semiVar} we made use of the stationarity assumption.
    
    In the next step we generate the empirical semivariogram based on the test vectors $V$ on $\mathcal{V}$ and the coordinate distance function $d^c(i,j)$. To this purpose, we collect all pairs of values $(d^c(i,j),(v_i^{(\ell)}-v_j^{(\ell)})^2), i,j\in \mathcal{V}$ and $\ell=1,\ldots,K$ in the so-called variogram cloud. Discretising the range of all values $d^c(i,j)$ into bins of width $\Delta$, we obtain the empirical semivariogram as
    \begin{equation}
    \label{eq:empSemivar}
        \hat\gamma (h)=\frac{1}{K\# \{(i,j):|h|-\frac{\Delta}{2}\leq d^c(i,j)<|h|+\frac{\Delta}{2}\}}\sum_{\{(i,j):|h|-\frac{\Delta}{2}\leq d^c(i,j)<|h|+\frac{\Delta}{2}\}\atop \ell=1,\ldots,K}(v_i^{(\ell)}-v_j^{(\ell)})^2.
    \end{equation}
    In the following, parametric semivariogram functions $\gamma_\theta$ are fitted to $\hat\gamma(|h|)$, mostly using weighted least squares \cite{bivand2008applied}. There are many known families of semivariograms, but in this work we use the exponential family 
    \begin{equation}\label{eq:spherical_model}
        \gamma_\theta(|h|)=\sigma^2\left(1-e^{-\left(\frac{|h|}{\eta}\right)}\right),~~\theta=(\sigma^2,\eta)\in\mathbb{R}_+^2
    \end{equation}
    and the spherical family defined by
    \begin{equation}\label{eq:exponential_model}
        \gamma_\theta(|h|)=\left\{\begin{array}{cc}\sigma^2\left(\frac{3|h|}{2a}-\frac{1}{2}\left(\frac{h}{\eta}\right)^3\right)&\mbox{for } |h|<\eta\\\sigma^2&\mbox{for }|h|\geq \eta\end{array}\right. , ~~\theta=(\sigma^2,\eta)\in\mathbb{R}_+^2.
    \end{equation}
    Empirical semivariograms and fitted semivariograms of exponential and spherical type can be seen in Figures \ref{fig:semivarCircle} and  \ref{fig:semivarSquare} below. For further models we refer to \cite{bivand2008applied,sherman2011spatial}. From the fitted semivariogram one then computes the covariance function $C_\theta(|h|)$ which can be used to compute $C_\mathcal{H}$. 
    
    Provided that the correlation length $\eta$ (also called range in the geostatistical literature) is much smaller than the size of the underlying domain, it is often enough to work with just a few or even just one test vector, $K=1$: If spatial correlations quickly decrease, the random field effectively contains many resamplings of its statistics in just one realization, i.e., in one test vector over a sufficiently extended grid $\mathcal{V}$. 
    
    \subsection{\label{sec:graphDist} Parametric Semivariograms for Inhomogeneity and Anisotropy}
   In most situations where homogeneity and isotropy of $X$ cannot be expected, we replace the coordinate metric $d^c(i,j)$ with the graph metric $d^A(i,j)$ introduced in subsection \ref{sec:locKrig}. This is very much in the spirit of algebraic multigrid approaches, where the knowledge of coordinate lists of the variables cannot be guaranteed and the use of the graph distance dates back at least to the definition of strength-of-connection in classical AMG~\cite{RugeStue1987}.  
   
   This pragmatic approach however comes with a conceptional problem. It is not clear if the weighted graph obtained from $\mathcal{V}$ and $A$ with distance function $d^A(i,j)$ can be isometrically embedded to some space $\mathbb{R}^{d'}$. In fact, this is possible if and only if the condition
   \begin{equation}
       \label{eq:embedding}
       w^T[d^A(i,j)^2]_{i,j\in\mathcal{V}}w\leq 0, \text{\ for all\ } w\in \mathbb{R}^n \mbox{ s.t. }\mathbbm{1}_n^Tw=0
   \end{equation}
   holds, see \cite{deza2009geometry}. If this is true, then the positivity of the covariance matrix $C$ follows from the fact that for many families of parametric semivariograms $C_\theta(|h|)$ defines a positive definite function in \emph{any} dimension $d$ and in particular in the dimension $d'$ of the isometric embedding, which can differ from the dimension $d$ of the underlying continuum problem.
   
   If the condition \eqref{eq:embedding} is violated, formula \eqref{eq:kriging} still defines an interpolator, as is easily checked, but the Kriging variance \eqref{eq:krigVar} needs no longer to be non-negative and the probabilistic interpretation of the Gaussian process in gone. In fact, we observe this in the examples presented below, despite observing high correlation between coordinate distance and graph distance.
   
   In this situation we still can carry out the algebraic manipulations from both formulae and we use the Gaussian process in the probabilistic sense simply as a source of inspiration. 
   
   Note however that local models of Gaussian processes on $\mathcal{C}_i\cup\{i\}$, $i\in \mathcal{F}$ with metric $d^A(i,j)$ may well exist as the embedding problem for such smaller graphs is much alleviated. In fact, in the numerical examples given below we do not observe any non-positive covariance metrics $\widehat C_{\mathcal{C}_i\cup\{i\}}$ for moderate size of $q_{\rm max}$. We therefore suggest that the calculation of local Kriging predictors and local Kriging variances can still be used and still give reasonable results, even though the (global) probabilistic interpretation has to be used with caution and strictly speaking the localized version of \eqref{eq:embedding} should be checked, at least if the method proposed in the following section does not show the expected performance. If this condition is violated, it seems to be reasonable to lower $q_{\rm max}$ at least locally, to obtain smaller local graphs which are more easily embedded. 
   
    \section{Adaptive Coarsening using Gaussian Processes}\label{sec:coarsening}
    The fact that algebraically smooth error can be interpreted as instances of a spatial Gaussian process allows us to use the methodology of Gaussian processes and the Kriging interpolation developed in~\cref{sec:gaussian} to solve the coarsening problem of algebraic multigrid methods. The calculation of a splitting of the variables $\mathcal{V}$ into a set of coarse variables $\mathcal{C}$ and the remaining fine variables $\mathcal{F} = \mathcal{V}\setminus \mathcal{C}$ as well as the computation of interpolation suited to this setting is carried out in three consecutive steps.
    
    First, we need to determine the covariance structure $C$ of algebraically smooth error when viewed as instances of a spatial Gaussian process. In order to do so, we start with a number, $K$, of test vectors $v^{(1)},\ldots,v^{(K)}$, each with entries that are normally distributed with mean zero and variance $1$. These initially random vectors are then subject to a number, $\nu$ of smoothing iterations with right-hand-side zero, where we employ the smoother that is going to be used in the algebraic multigrid method~\eqref{eq:eprop} as well. These test vectors are then fed into the calculation of the covariance structure. In our tests we compare localized non-parametric approaches (cf.~\cref{sec:covMod}), which use a sizeable number of test vectors, with parametric models (cf.~\cref{sec:semiVar}) that make use of only a small number of test vectors.
    
    Once the covariance structure is determined we can use Kriging interpolation to tackle the coarsening problem. The set of coarse variables $\mathcal{C}$ has to be chosen in such a way that interpolation of information from these variables to the remaining variables $\mathcal{F}$ is as accurate as possible for algebraically smooth error. In accordance with the interpretation of Kriging interpolation as the interpolation that minimizes the MSE under the assumption that algebraically smooth error can be interpreted as instances of a Gaussian field, we can use the variance of the Kriging estimator~\eqref{eq:krigVar} in order to define the coarse variable set. Starting with $\mathcal{C} = \emptyset$ and using the fact that any variable that ends up in $\mathcal{C}$ during the coarsening has zero variance after interpolation we proceed to add those variables to $\mathcal{C}$ with largest variance. In case there is a tie, we choose the first occurrence, but other selection strategies are possible as well. 
    Exploiting the fact that the correlation distance can be used to limit the reach of the Kriging interpolation, we can actually add multiple variables to $\mathcal{C}$ at the same time if they are spaced so far apart that interpolation between them is not considered due to the localization of the Kriging interpolation. After variables are added to $\mathcal{C}$, we update the Kriging interpolation of all affected variables. By using one of the pseudo-distances we first determine for each $i\in \mathcal{F}$ the set of interpolatory variables $\mathcal{C}_i \subseteq \mathcal{C}$, again respecting the limitation of reach due to a finite correlation distance. Based on these sets, the Kriging estimator and the corresponding variances or MSE are calculated/updated. This process of adding variables to $\mathcal{C}$ based on the uncertainty with which we can predict the value at the respecting variable and updating Kriging interpolation is repeated until a prescribed tolerance on either the size of the coarse variable set $\mathcal{C}$, typically relative to the total number of variables $\mathcal{V}$, or the largest remaining uncertainty of the Kriging estimator is reached. The resulting process is roughly summarized in~\cref{alg:coarsening}. 
    
    \begin{algorithm2e}[ht]
    \SetAlgorithmStyle
    \caption{Coarsening based on Kriging interpolation}\label{alg:coarsening}
    \KwData{}\medskip
        Initialize $\mathcal{C} = \emptyset$, $\mathcal{F} = \mathcal{V}$\;
        \While{$|\mathcal{C}| < n_c$}{
            Choose $i \in \mathcal{F}$ with largest variance $\mu_i$\;
            Add $i$ to $\mathcal{C}$, remove $i$ from $\mathcal{F}$\;
            \For{$j\in \mathcal{F}$ with $i \in \mathcal{C}_j$}{
                Update Kriging interpolation with new $\mathcal{C}_j$ set\;
                Compute corresponding updated variances $\mu_j$ (set $\mu_i = 0$)\;
            }
        }
    \end{algorithm2e}
    
    \section{Numerical Case Studies}\label{sec:numerics}
    

    In order to gauge the efficiency of the new coarsening scheme based on Gaussian fields and Kriging interpolation we consider the general diffusion problem
    \begin{equation*}\label{eq:poisson_univ}  
    - \left( c_1 \frac{\partial^2}{\partial x^2} + c_2 \frac{\partial^2}{\partial y^2} + c_3 \left( \frac{\partial}{\partial x} \frac{\partial}{\partial y} + \frac{\partial}{\partial y}  \frac{\partial}{\partial x} \right) \right) u = f \, ,
    \end{equation*} for constant and anisotropic coefficients by choosing $c_1,c_2$ and $c_3$ accordingly. We further choose the computational domain as the unit square $(0,1)^2$, employing a finite difference discretization on a regular mesh and the unit circle $U(1)$, where the discretization is defined by linear finite elements on a triangularization generated in MATLAB. We considered $4$ parameter combinations in our tests collected in~\cref{tab:parameters}.
    \begin{table}[ht]\small
        \centering
        \begin{tabular}{c||c|c|c|c}
            name & s-iso & s-aniso & c-iso & c-aniso \\\hline
            domain & square & square & circle & circle \\
            \textbf{$n$} & $2025$ & $2025$ & $2521$ & $2521$\\
            $\begin{bmatrix} c_1 & c_3 \\ c_3 & c_2\end{bmatrix}$ & $\begin{bmatrix} 1 & 0 \\ 0 & 1\end{bmatrix}$ & $\begin{bmatrix} 1 & 0 \\ 0 & 10^{-2}\end{bmatrix}$ & $\begin{bmatrix} 1 & 0 \\ 0 & 1\end{bmatrix}$ & $\begin{bmatrix} 1 & 0 \\ 0 & 10^{-2}\end{bmatrix}$
        \end{tabular}
        \caption{Parameter choices of the considered test cases.}
        \label{tab:parameters}
    \end{table}
    
    The overall aim of these case studies is to get a first impression of the performance of the new approach of constructing coarse variable sets and interpolation. Thus we compare methods with respect to the underlying covariance model, i.e., the empirical covariance function vs.~semivariogram models (spherical and exponential). In addition, we test the robustness of the approach with respect to the number of test vectors used to construct the covariance model. In all tests we employ localization in the calculation of Kriging interpolation, which is an inevitable technique to guarantee overall linear complexity and well-posedness of the Kriging interpolation as discussed in~\cref{sec:locKrig} especially for the empirical construction.
    
    The test vectors are generated by applying just one iteration of a colored Gauss-Seidel iteration to white noise vectors and we use the same iterative method for the smoother of our resulting two-grid method as well. That is, in the reported results we run a $V(1,1)$-cycle two-grid method with a direct solve for the coarse grid system of equations.
    
    \paragraph{Pseudo-distances}
    The use of pseudo-distances is an important aspect when it comes to the independence from problem specific knowledge as is required by a true algebraic multigrid approach. Thus we first compare the correlation of the graph pseudo-distance and the true geometric distances based on the coordinates of the unknowns. In this we found a correlation of the distances for the cases c-iso of $97.68\%$ and for s-iso of $99.65\%$, which leads us to believe that using the algebraic graph distance as a pseudo distance in our Gaussian field analysis is viable.

\begin{figure}[h]
    \centerline{
    \includegraphics[width=.45\textwidth]{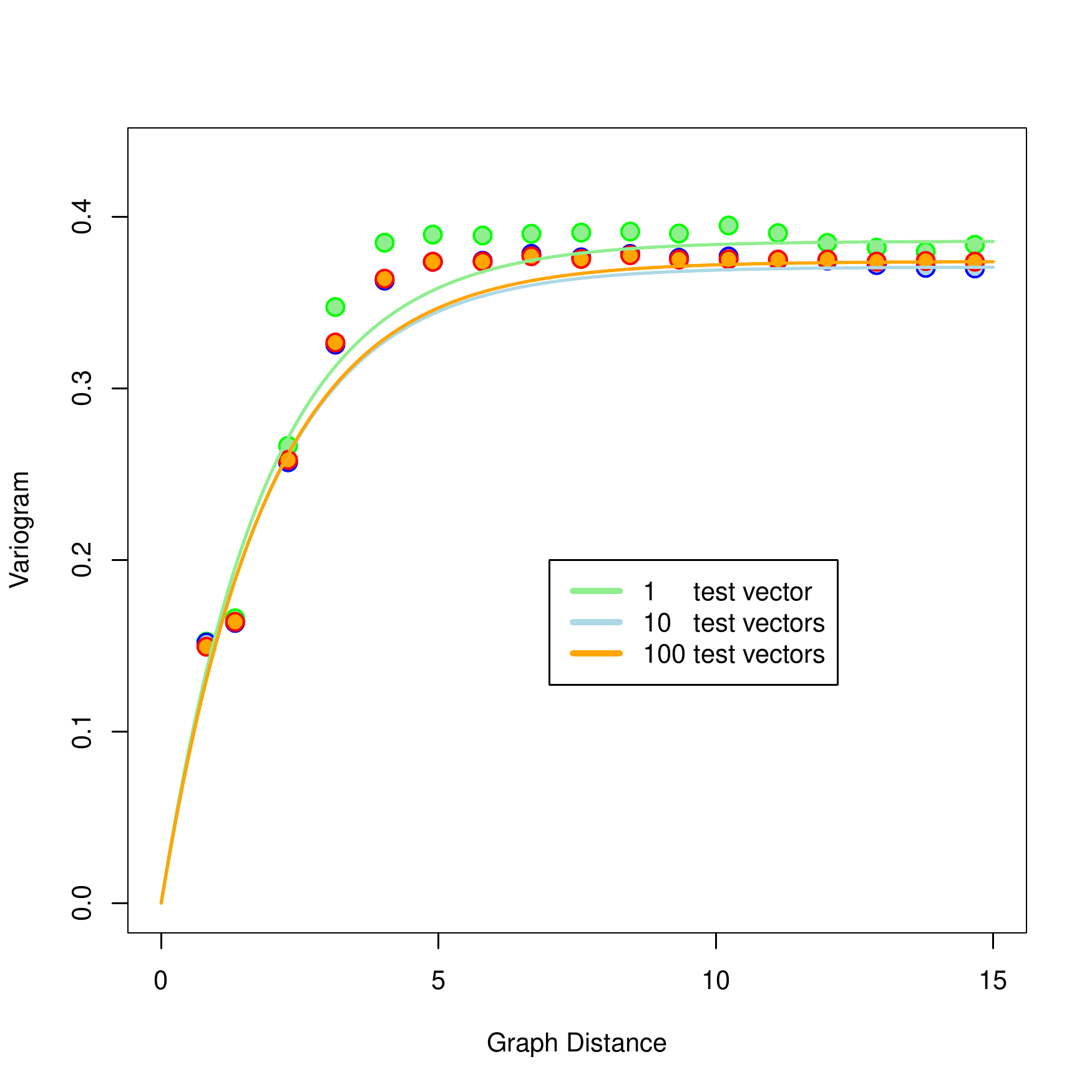}\hfill
    \includegraphics[width=.45\textwidth]{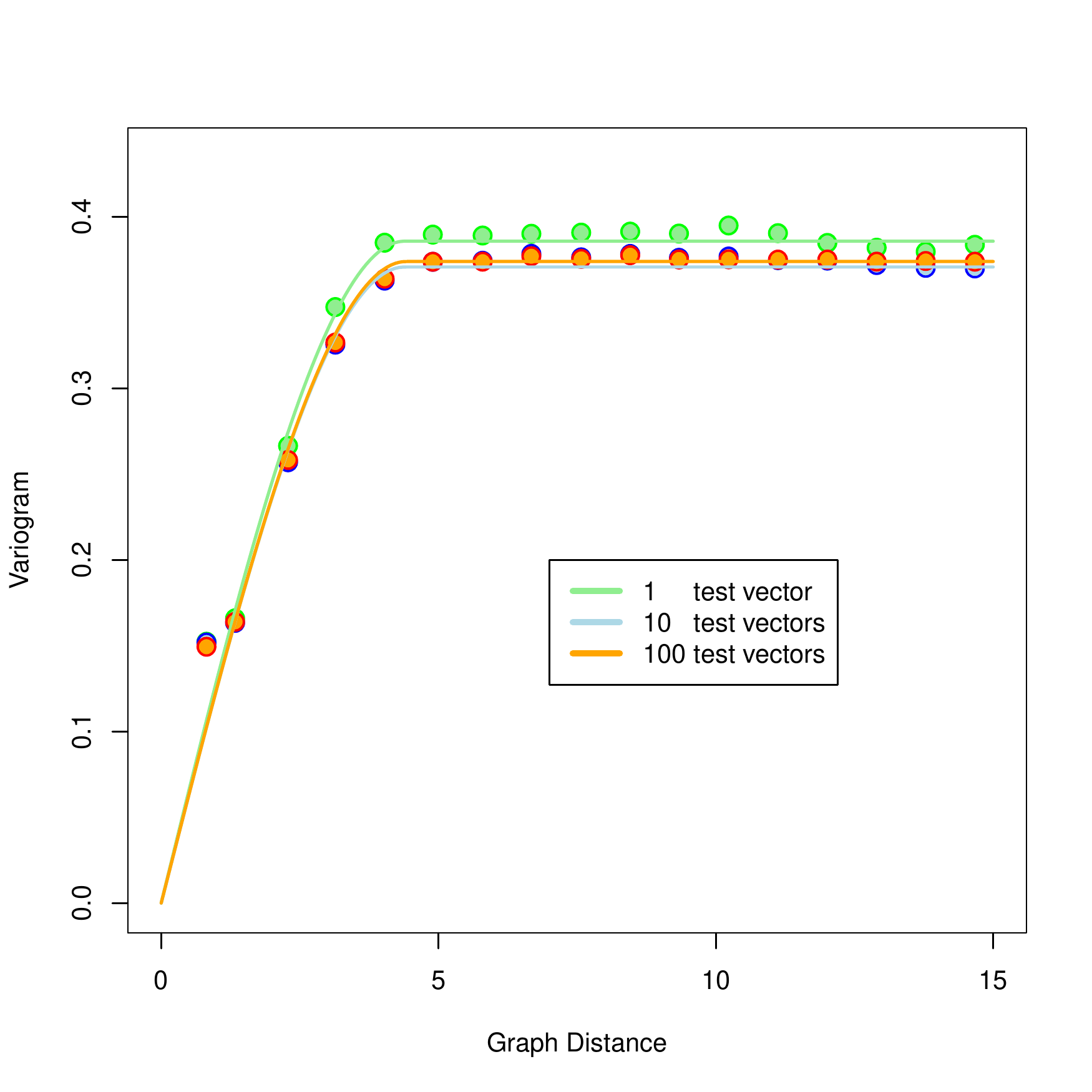}
    }
    \centerline{
    \includegraphics[width=.45\textwidth]{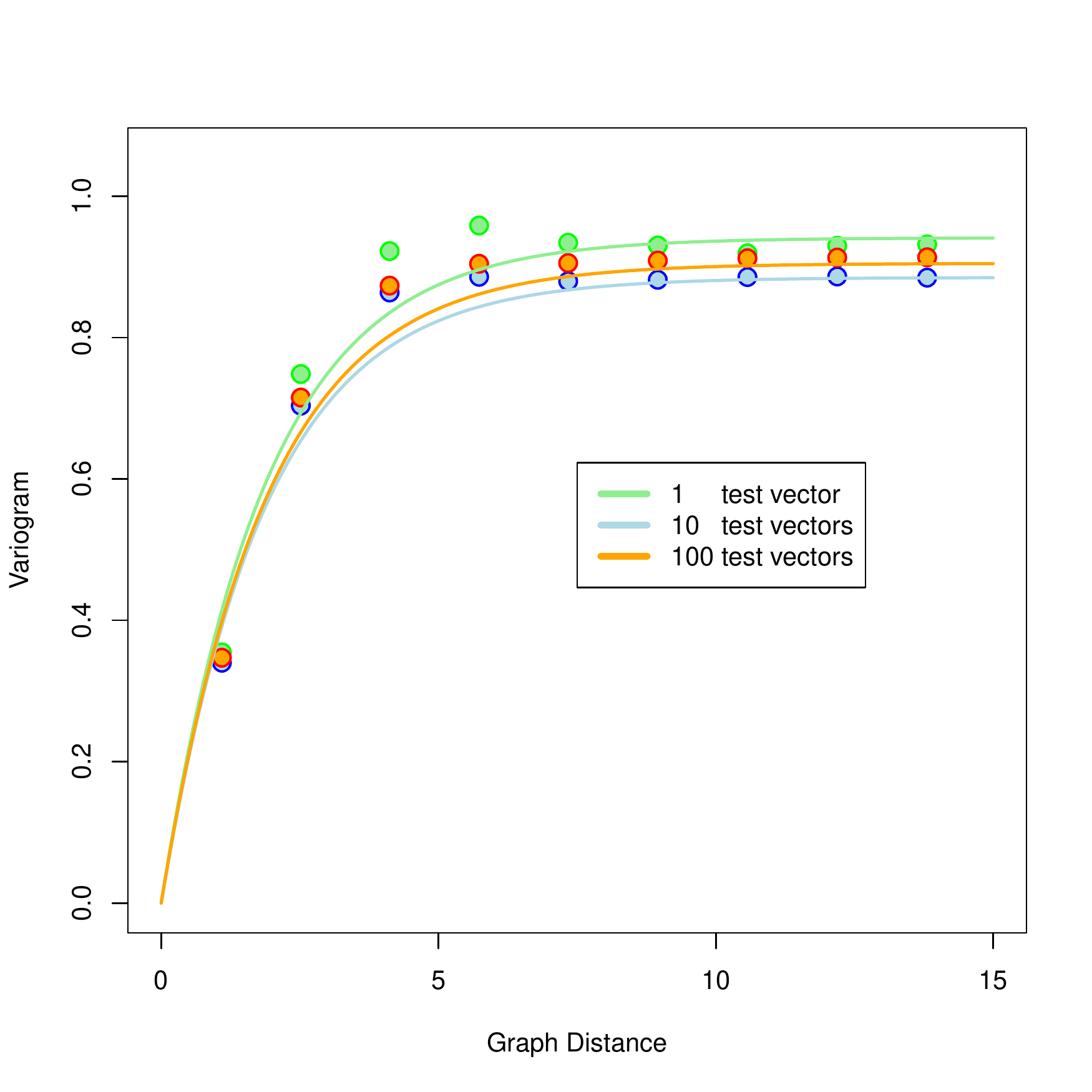}\hfill
    \includegraphics[width=.45\textwidth]{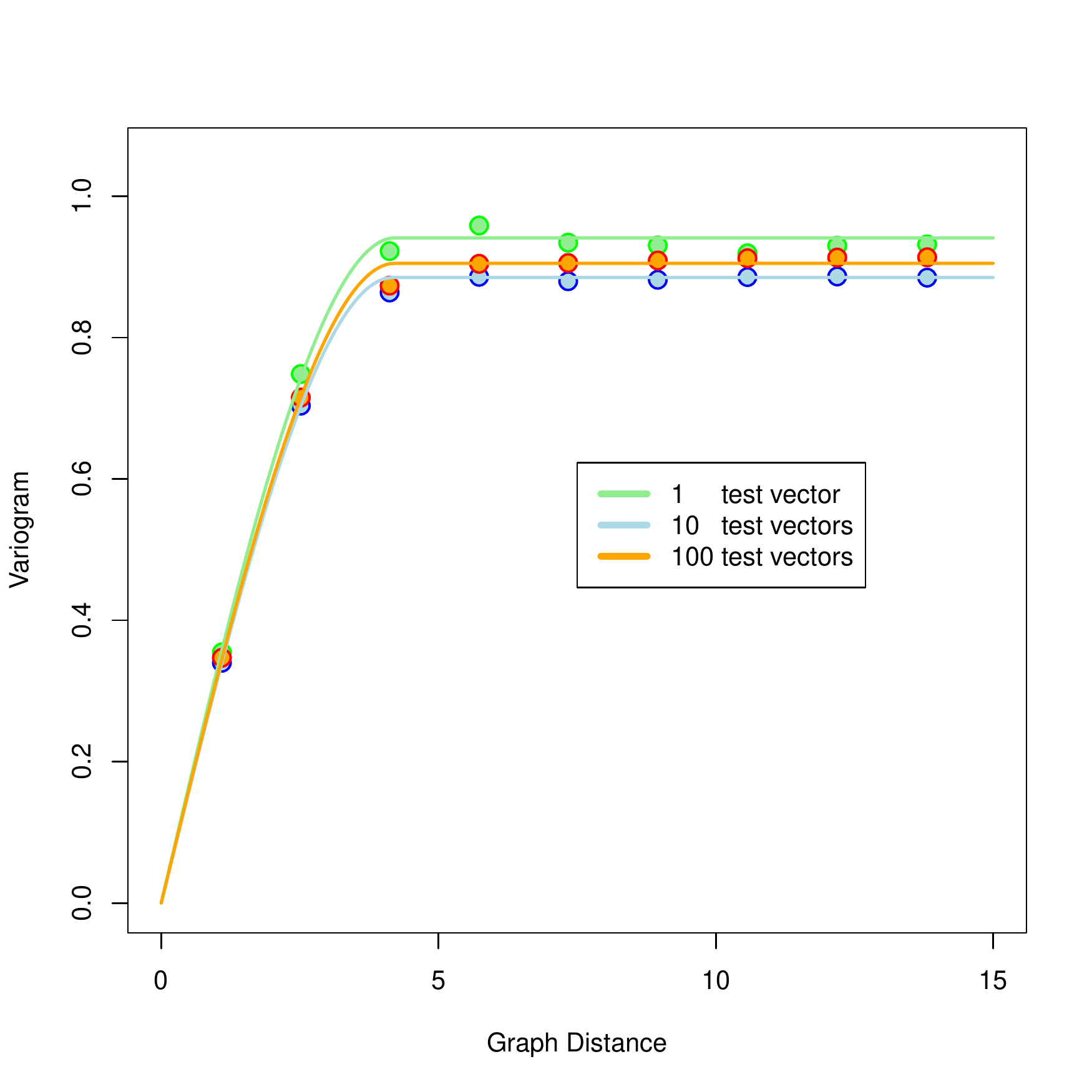}
    }
    \caption{\label{fig:semivarCircle} Fitted semivariogram models for the unstructured mesh on the circle: The top row displays the isotropic  and the bottom row the anisotropic case. On the left the exponential semivariogram model is used and on the right the spherical.}
    \end{figure}

\begin{figure}[h]
    \centerline{
    \includegraphics[width=.45\textwidth]{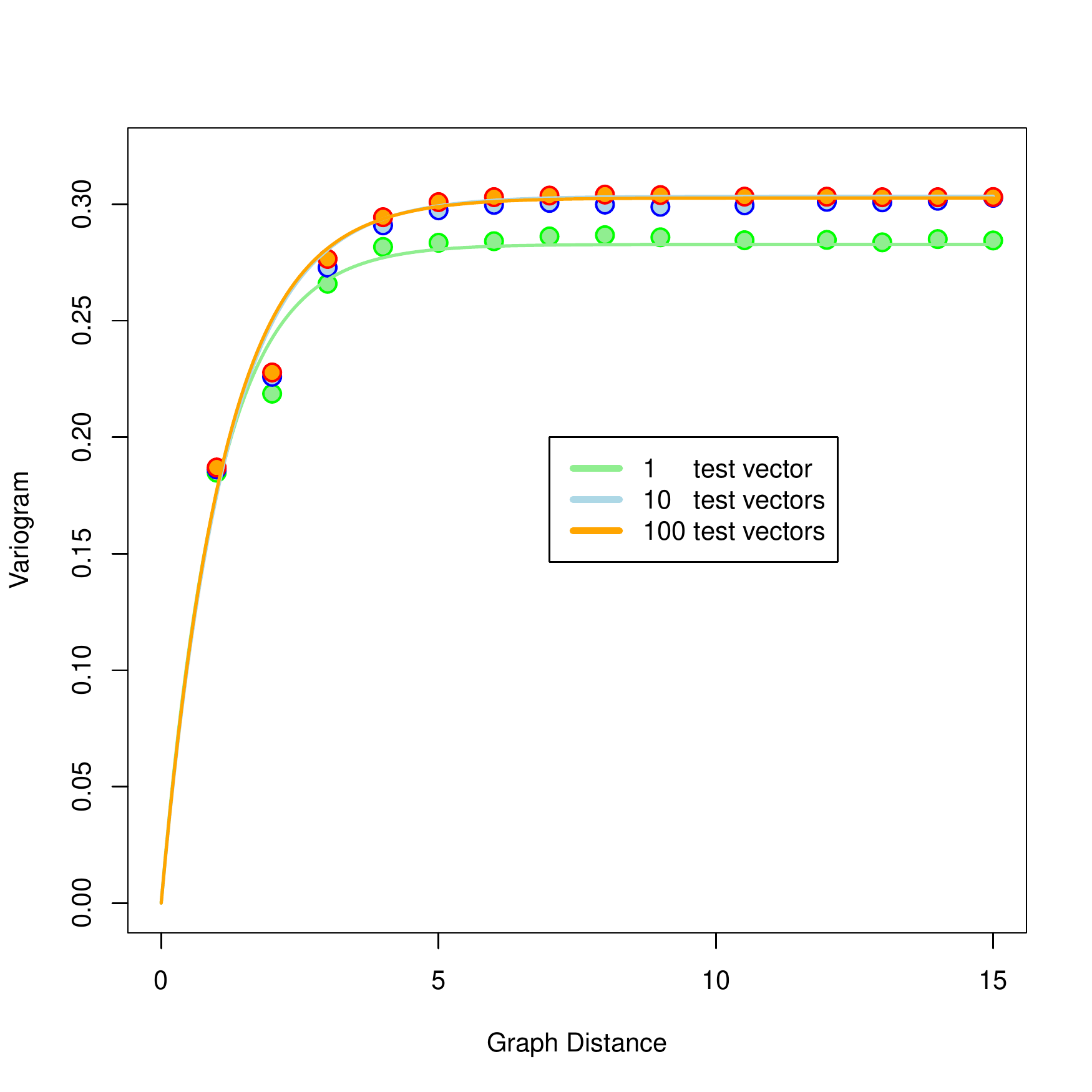}\hfill
    \includegraphics[width=.45\textwidth]{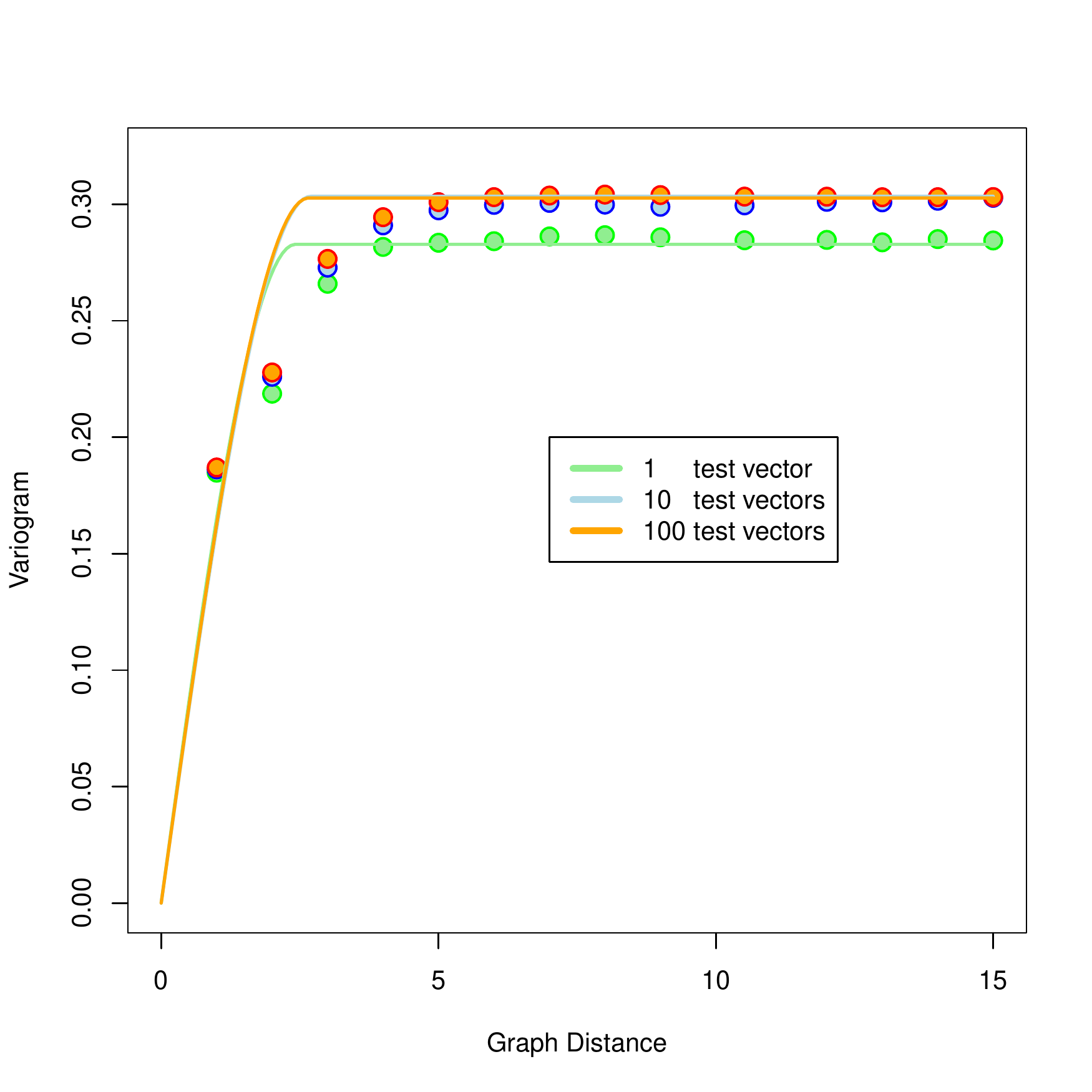}
    }
    \centerline{
    \includegraphics[width=.45\textwidth]{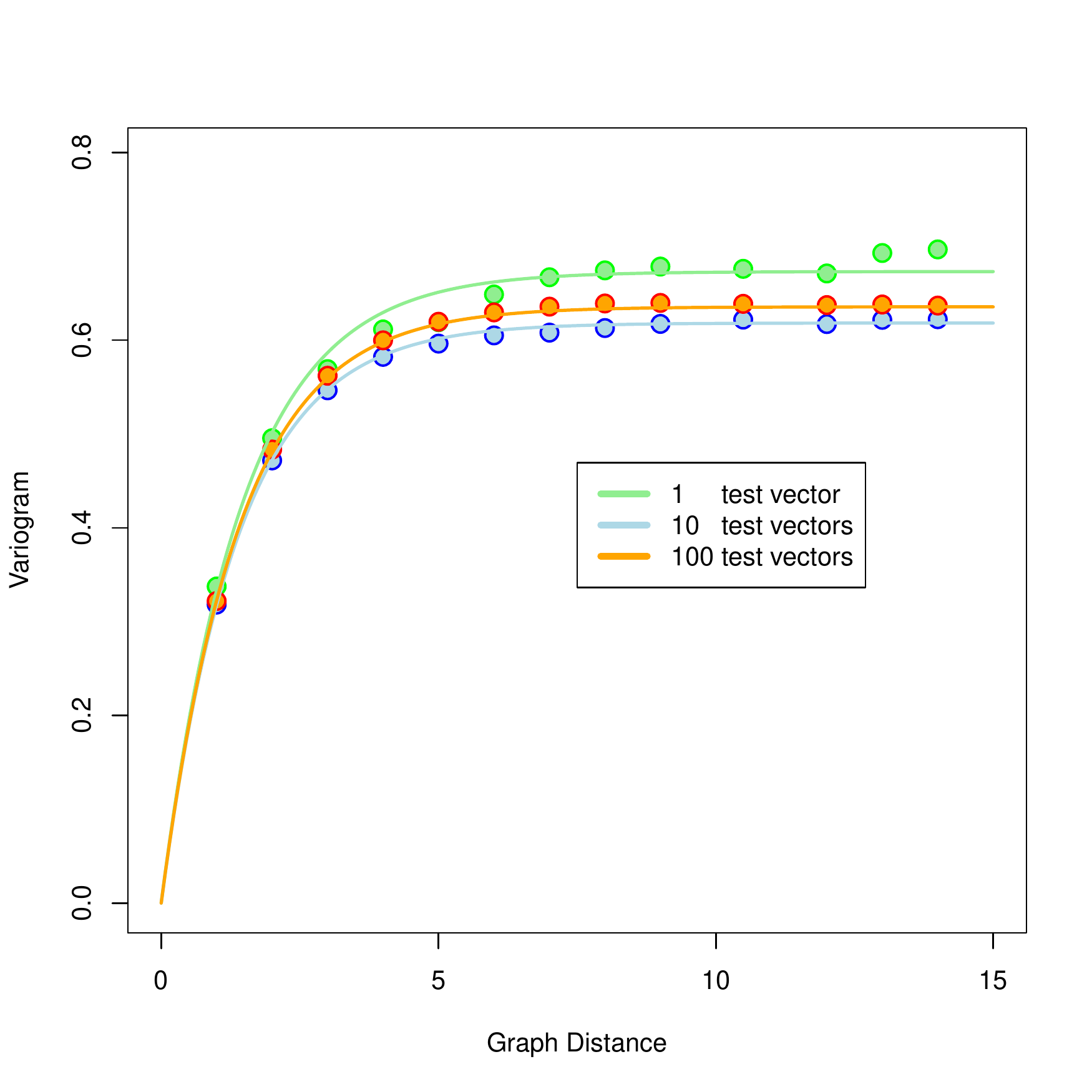}\hfill
    \includegraphics[width=.45\textwidth]{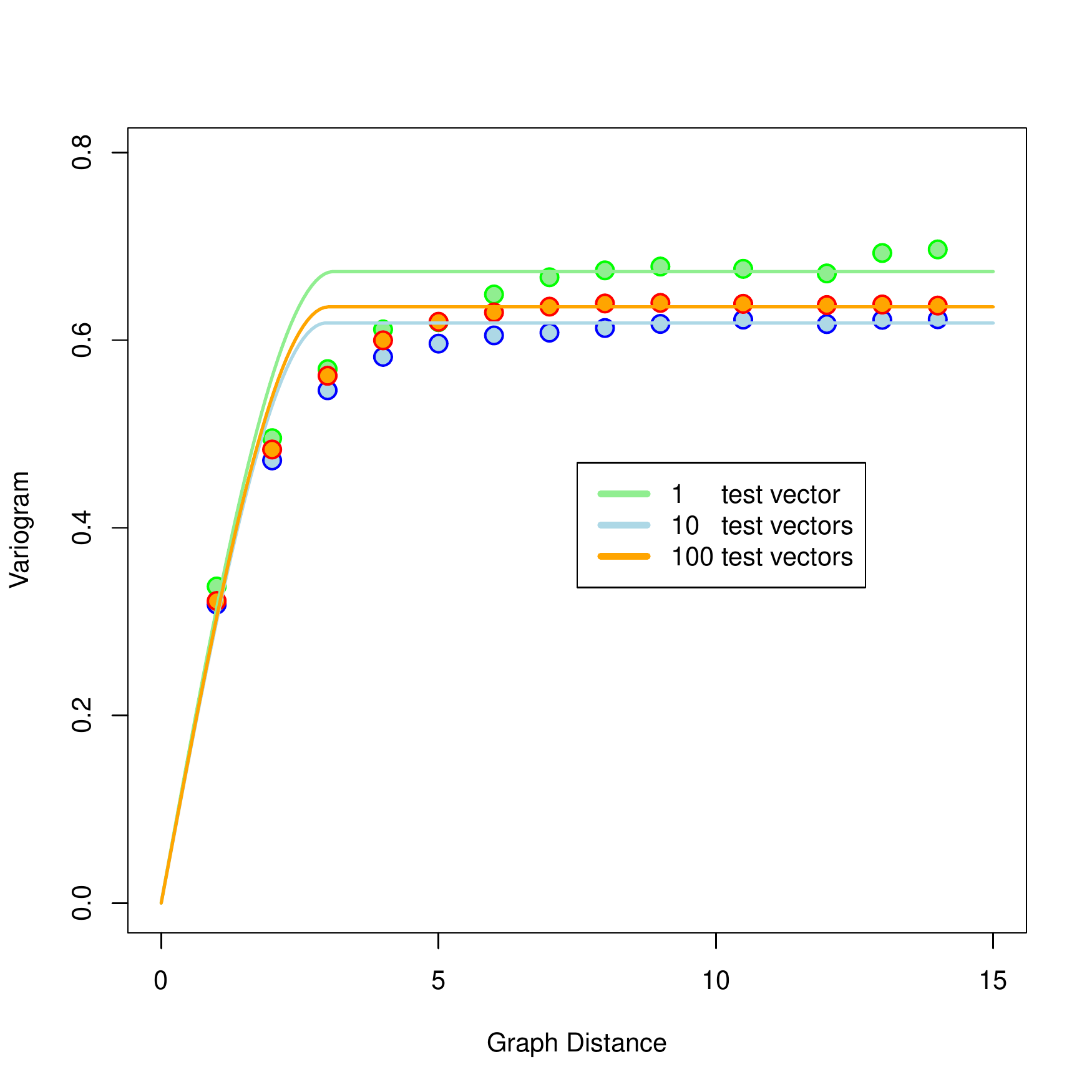}
    }
    \caption{\label{fig:semivarSquare} Fitted semivariogram models for the structured mesh on the square: The top row displays the isotropic  and the bottom row the anisotropic case. On the left the exponential semivariogram model is used and on the right the spherical.}
    \end{figure}

    
    \paragraph{Covariance models}
    In the following we fit the exponential  and the shperical  semivariogram models to the empirical semivariogram  generated from 1, 10 and 100 test vectors, see equations \eqref{eq:exponential_model}, \eqref{eq:spherical_model} and \eqref{eq:empSemivar}. This is done for the isotropic and the anisotropic case, both for the circle (Figure \ref{fig:semivarCircle}) and the square grid (Figure \ref{fig:semivarSquare}). The fits and models are performed using the R library \texttt{gstat}. The fits expose reasonable quality and not much variation caused by the number of test vectors used.
    
    \paragraph{Coarsening and two-grid results}
    Finally, after analyzing the components of the metrics and fits underlying the Gaussian process framing of the coarse grid construction we are ready to apply our coarsening algorithm and present two-grid results using the Kriging interpolation. In this we use the test cases described in~\cref{tab:parameters} combined with 
    \begin{itemize}
        \item an empirical construction of the covariance structure, which we term \emph{emp}-$K$, 
        \item a spherical covariance model based on a semivariogram fit~\eqref{eq:spherical_model}, termed \emph{sph}-$K$,
        \item an exponential covariance model based on a semivariogram fit~\eqref{eq:exponential_model}, termed \emph{exp}-$K$.
    \end{itemize} In this naming convention $K$ represents the number of test vectors used in the construction of the model or semivariogram, respectively.

\begin{figure}[ht]
        \begin{center}\resizebox{.8\textwidth}{!}{
            \begin{tikzpicture}
                \begin{scope}
                    \node at (0,0) {\includegraphics[width=.45\textwidth]{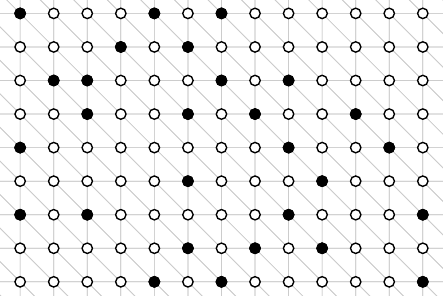}};
                \end{scope}
                
                \begin{scope}[xshift = .5\textwidth]
                    \node at (0,0) {\includegraphics[width=.45\textwidth]{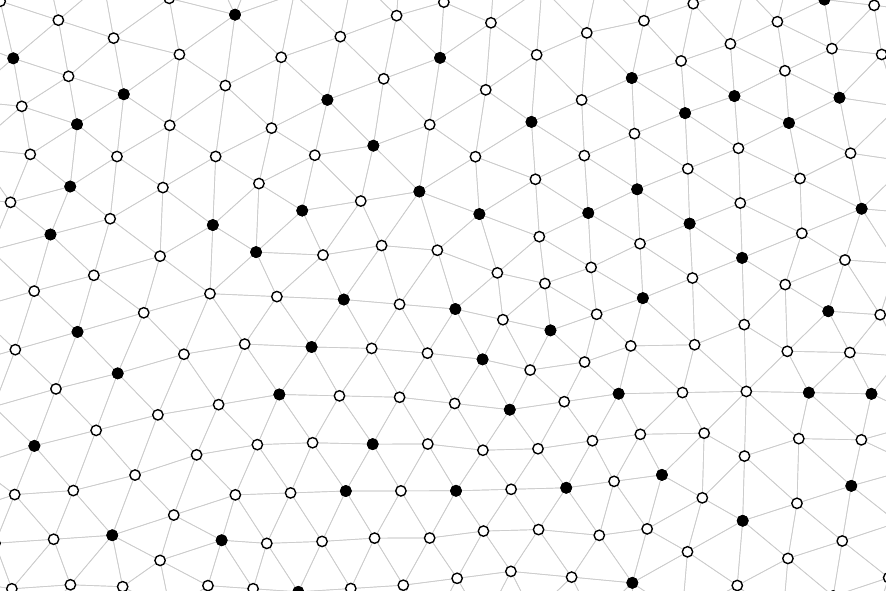}};
                \end{scope}
                
                \begin{scope}[yshift = -.35\textwidth]
                    \node at (0,0) {\includegraphics[width=.45\textwidth]{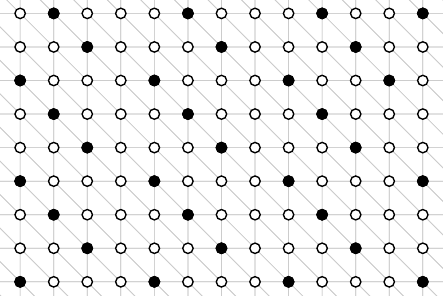}};
                \end{scope}
                
                \begin{scope}[xshift=.5\textwidth,yshift=-.35\textwidth]
                    \node at (0,0) {\includegraphics[width=.45\textwidth]{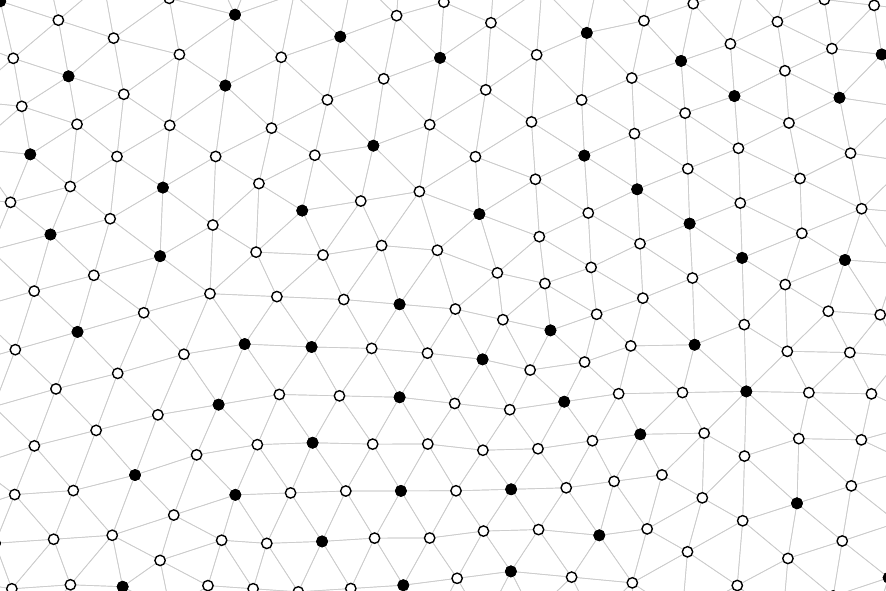}};
                \end{scope}
                
                \begin{scope}[yshift=-.7\textwidth]
                    \node at (0,0) {\includegraphics[width=.45\textwidth]{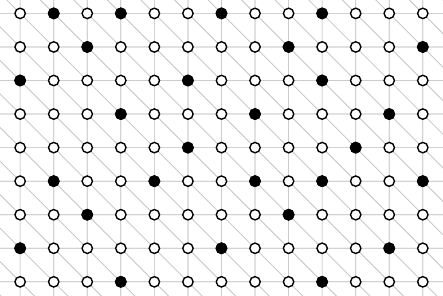}};
                \end{scope}
                
                \begin{scope}[xshift=.5\textwidth,yshift=-.7\textwidth]
                    \node at (0,0) {\includegraphics[width=.45\textwidth]{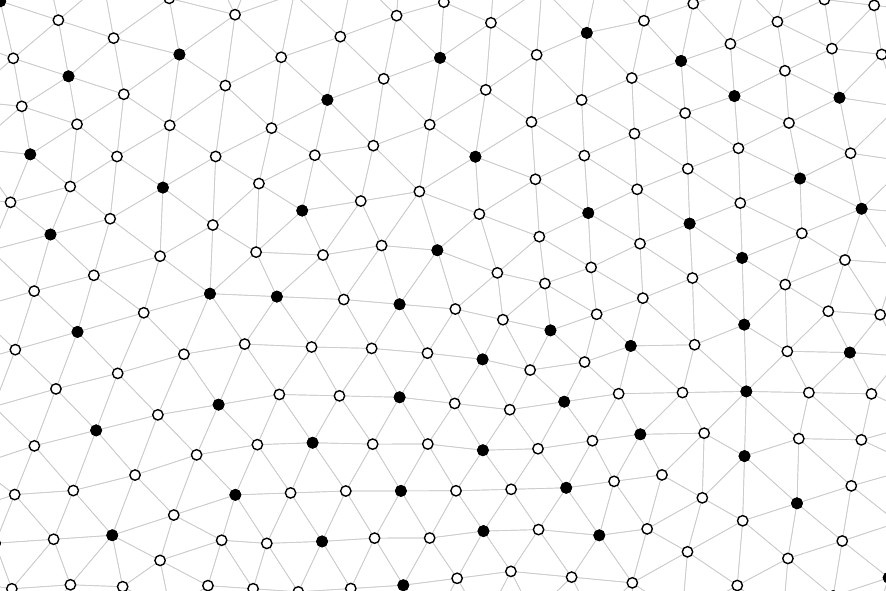}};
                \end{scope}
            \end{tikzpicture}}
        \end{center}
        \caption{Illustration of the variable splittings obtained by emp-$10$ (top), sph-$1$ (mid) and exp-$1$ (bottom) for the isotropic test cases on the square (left) and circle (right). In order to avoid cluttering of the illustrations we just show a section of the bulk of the domains (i.e., $[.35,.4]\times[.65,.6]$ and $[-.3,-.2]\times[.3,.2]$ for square and circle, resp.).}\label{fig:corsening-iso}
    \end{figure}

    In~\cref{tab:twogrid-iso} we collect results of the resulting two-grid methods for the isotropic test cases s-iso and c-iso. The quantities we report are the asymptotic convergence rates, $\rho$ of the two-grid methods, which are indicative of the overall compatibility of the smoother and coarse-grid correction, and iteration numbers, $k$, of the conjugate gradients method preconditioned with the two-grid method that are required to reduce the initial residual by a factor of $10^8$. This latter quantity gives insight into the flaws of the two-grid construction. Oftentimes, as also explored in~\cite{BranBranKahlLivs2011,BranBranKahlLivs2015a,KahlRott2018}, the two-grid method might show bad asymptotic convergence rates even though the preconditioned conjugate gradients iteration converges rapidly. This typically corresponds to the presence of a few outliers in the spectrum of the preconditioned matrix and indicates that the two-grid construction provides better complementarity of coarse grid construction and smoothing than the asymptotic convergence rate suggests.
    \begin{table}[ht]
        \resizebox{\textwidth}{!}{
        \begin{tabular}{c}
        \begin{tabular}{c|c|c|c|c|c|c|c|c}
            \textbf{s-iso} & emp-$10$ & emp-$100$ & sph-$1$ & sph-$10$ & sph-$100$ & exp-$1$ & exp-$10$ & exp-$100$ \\\hline
             $\rho$ & $.387$ & $.302$ & $.256$ & $.251$ & $.253$ & $.224$ & $.225$ & $.222$\\
             $k$ & $10$ & $9$ & $9$ & $9$ & $10$ & $9$ & $8$ & $9$\\
        \end{tabular}\\ \\
        \begin{tabular}{c|c|c|c|c|c|c|c|c}
             \textbf{c-iso} & emp-$10$ & emp-$100$ & sph-$1$ & sph-$10$ & sph-$100$ & exp-$1$ & exp-$10$ & exp-$100$ \\\hline
             $\rho$ & $.563$ & $.319$ & $.275$ & $.27$ & $.273$ & $.314$ & $.294$ & $.303$\\
             $k$ & $14$ & $10$ & $10$ & $10$ & $10$ & $11$ & $10$ & $10$\\
        \end{tabular}
        \end{tabular}
        }
        \caption{Asymptotic convergence rates $\rho$ of the two-grid V($1,1$) cycle and number of iterations $k$ of the conjugate gradients method preconditioned with the two-grid method to reduce the initial residual by a factor of $10^{8}$ for both isotropic test cases. All approaches generate a coarse variable set with $n_c = \tfrac{n}{4}$ variables, use a localization of distance $4$ and a caliber of $4$, i.e., $|\mathcal{C}_i|\leq 4$ for all $i\in \mathcal{F}$.}
        \label{tab:twogrid-iso}
    \end{table}
    Taking into account that all results of~\cref{tab:twogrid-iso} are generated with test vectors that are smoothed by only a single iteration of colored Gauss-Seidel and that a $V(1,1)$-cycle is employed the results are surprisingly good; cf.~\cite{BranBranKahlLivs2011}. In part this can be explained by the implicit preservation of the constant vector in the Kriging interpolation, a modification that has been shown to be particularly effective in improving the performance of bootstrap AMG in~\cite{BranBranKahlLivs2011,BranBranKahlLivs2015b,MantMcCoParkRuge2010}. The most interesting observation is the fact that the model based constructions of Kriging interpolation are competitive when using only a single test vector. Due to the fact that the descriptive power of the semivariogram approach should get better the more points are available for sampling, i.e., the finer the discretization, this finding should scale well with respect to the problem size. In~\cref{fig:corsening-iso} we collected some sections of the resulting coarsenings of the methods considered in~\cref{tab:twogrid-iso}. 
    
    The coarsenings do not show any particularly interesting features. As only sections in the bulk are shown and the coarsening ratio of the depicted section is smaller than the preset $\tfrac{1}{4}$, a mild agglomeration of coarse grid variables at the boundary takes place. It remains to be seen if this poses a problem when recursing on the construction in a multigrid fashion.
    
\begin{figure}[ht]
        \begin{center}\resizebox{.8\textwidth}{!}{
            \begin{tikzpicture}
                \begin{scope}
                    \node at (0,0) {\includegraphics[width=.45\textwidth]{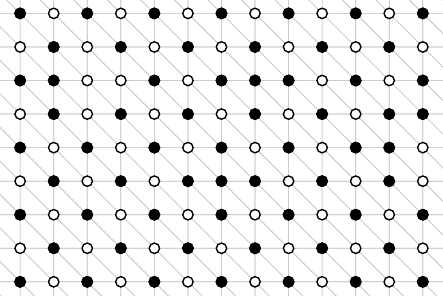}};
                \end{scope}
                
                \begin{scope}[xshift = .5\textwidth]
                    \node at (0,0) {\includegraphics[width=.45\textwidth]{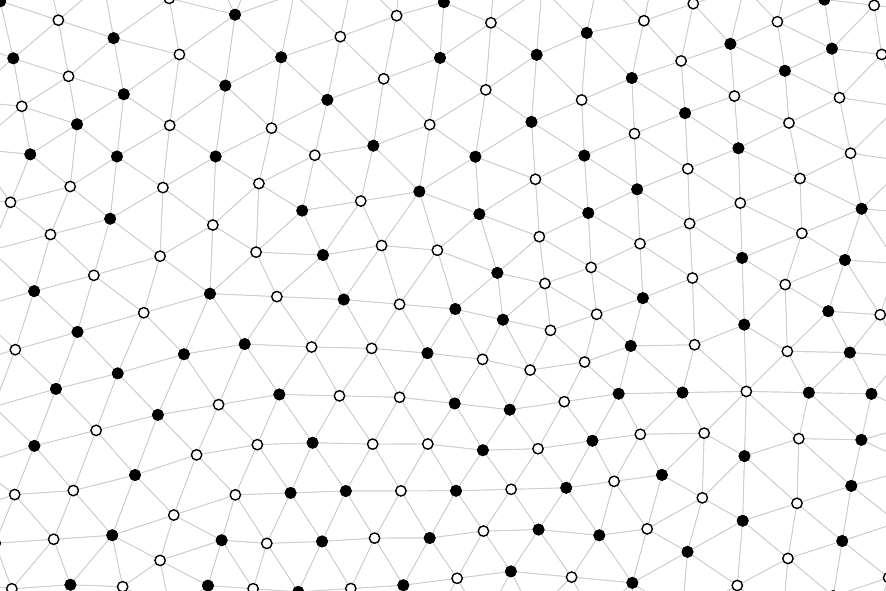}};
                \end{scope}
                
                \begin{scope}[yshift = -.35\textwidth]
                    \node at (0,0) {\includegraphics[width=.45\textwidth]{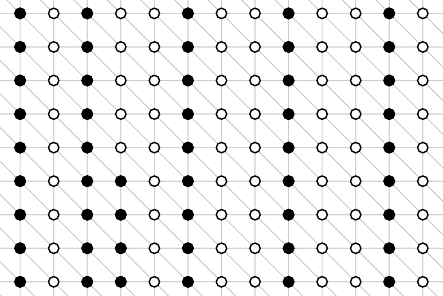}};
                \end{scope}
                
                \begin{scope}[xshift=.5\textwidth,yshift=-.35\textwidth]
                    \node at (0,0) {\includegraphics[width=.45\textwidth]{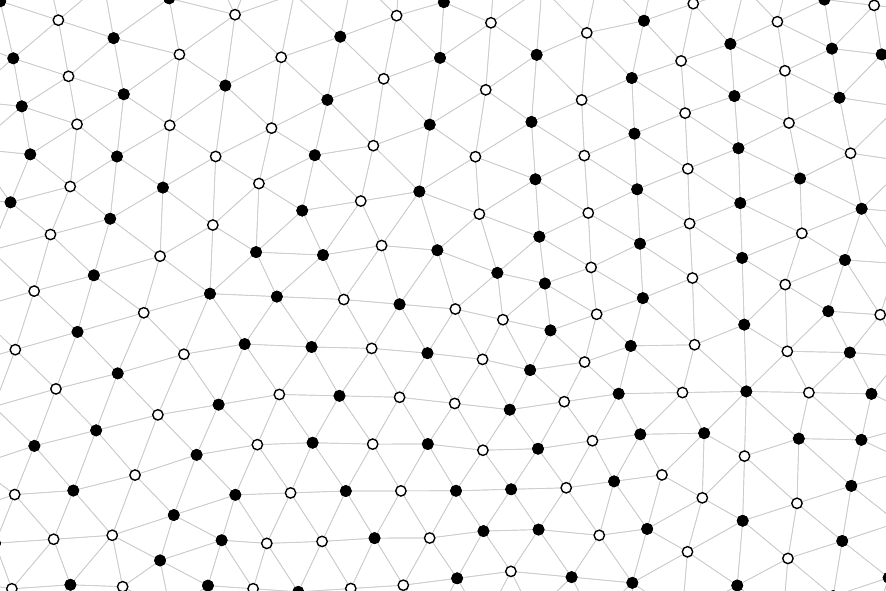}};
                \end{scope}
                
                \begin{scope}[yshift=-.7\textwidth]
                    \node at (0,0) {\includegraphics[width=.45\textwidth]{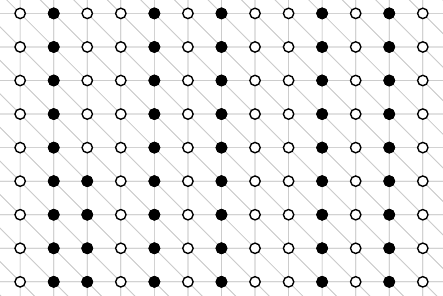}};
                \end{scope}
                
                \begin{scope}[xshift=.5\textwidth,yshift=-.7\textwidth]
                    \node at (0,0) {\includegraphics[width=.45\textwidth]{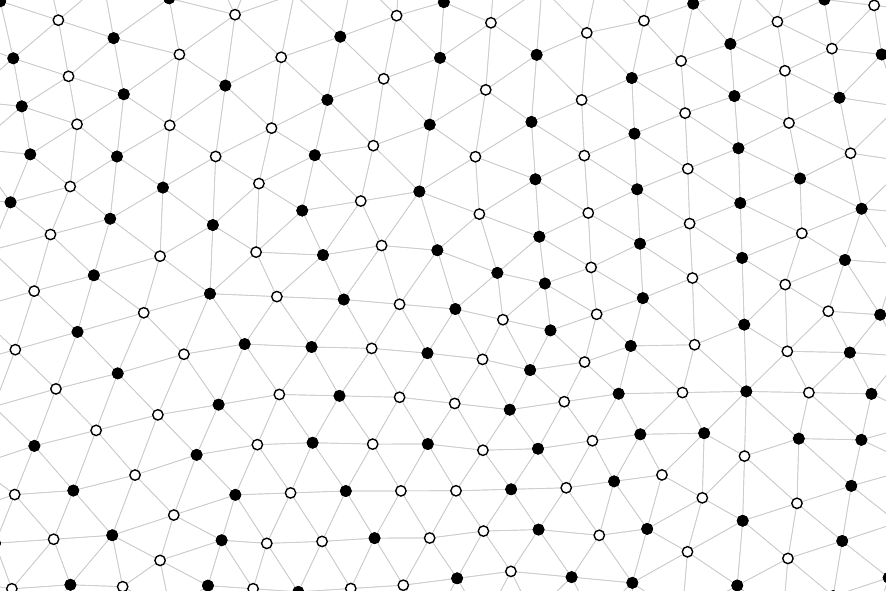}};
                \end{scope}
            \end{tikzpicture}}
        \end{center}
        \caption{Illustration of the variable splittings obtained by emp-$10$ (top), sph-$1$ (mid) and exp-$1$ (bottom) for the anisotropic test cases on the square (left) and circle (right). In order to avoid cluttering of the illustrations we just show a section of the bulk of the domains (i.e., $[.35,.4]\times[.65,.6]$ and $[-.3,-.2]\times[.3,.2]$ for square and circle, resp.).}\label{fig:corsening-aniso}
    \end{figure}

    The more interesting test cases arise when anisotropy is present in the model, especially for discretization on unstructured grids, where a canonical coarsening that follows the anisotropy is not available. Testing automatic coarsening approaches on the simple grid-aligned case in order to gauge their robustness and ability to reproduce the canonical coarsenings has been common in the past~\cite{BranBranKahlLivs2015b,KahlRott2018}. Analogous to the results for the isotropic test cases we report in~\cref{tab:twogrid-aniso} asymptotic convergence rates and iterations counts of the preconditioned conjugate gradients method. 
    \begin{table}[ht]
        \resizebox{\textwidth}{!}{
        \begin{tabular}{c}
        \begin{tabular}{c|c|c|c|c|c|c|c|c}
             \textbf{s-aniso} & emp-$10$ & emp-$100$ & sph-$1$ & sph-$10$ & sph-$100$ & exp-$1$ & exp-$10$ & exp-$100$ \\\hline
             $\rho$ & $.463$ & $.305$ & $.06$ & $.06$ & $.154$ & $.224$ & $.225$ & $.222$\\
             $k$ & $9$ & $8$ & $6$ & $6$ & $6$ & $9$ & $8$ & $9$\\
        \end{tabular}\\ \\
        \begin{tabular}{c|c|c|c|c|c|c|c|c}
             \textbf{c-aniso} & emp-$10$ & emp-$100$ & sph-$1$ & sph-$10$ & sph-$100$ & exp-$1$ & exp-$10$ & exp-$100$ \\\hline
             $\rho$ & $.648$ & $.533$ & $.69$ & $.684$ & $.681$ & $.704$ & $.71$ & $.685$\\
             $k$ & $19$ & $15$ & $22$ & $21$ & $21$ & $22$ & $24$ & $22$\\
        \end{tabular}
        \end{tabular}
        }
        \caption{Asymptotic convergence rates $\rho$ of the two-grid V($1,1$) cycle and number of iterations $k$ of the conjugate gradients method preconditioned with the two-grid method to reduce the initial residual by a factor of $10^{8}$. All approaches generate a coarse variable set with $n_c = \tfrac{n}{2}$ variables, use a localization of distance $4$ in graph distance  and a caliber of $2$ and $3$ for all $i\in \mathcal{F}$ for the test cases formulated on the square and circle, respectively.}
        \label{tab:twogrid-aniso}
    \end{table} The grid aligned anisotropy in the square test case s-aniso does not pose any difficulty for any of the model based approaches and again a single test vector is sufficient to obtain a good enough statistics for a suitable fit of the model as already suggested by~\cref{fig:semivarSquare,fig:semivarCircle}. The results suggest that the spherical model is better suited for this problem than the exponential model and we see a clear advantage of the model based approaches over the empirical approach, showing extremely fast convergence. Similar to the isotropic case we see a notable improvement when increasing the number of the test vectors for the empirical model from $10$ to $100$, which is quite frankly an unfeasible number of test vectors, but serves the illustrative purpose quite well. As can be seen in~\cref{fig:corsening-aniso} the coarsenings obtained by the different approaches for the square test case show that the anisotropy has been clearly detected and the coarsening constructed accordingly. When it comes to the results for the anisotropic problem on the circular domain and unstructured grid the results are comparable to results reported, e.g., in~\cite{BranBranKahlLivs2015b} for non-grid aligned anisotropies. Whilst all approaches yield good preconditioners for the conjugate gradients method with comparable iteration numbers, the asymptotic convergence rates are considerably worse at around $.7$ compared to the grid-aligned case for the model based approaches. Interestingly the empirical construction is able to cope with the unstructured grid better than the model based approaches. Taking into account that the model based approaches implicitly assume shift invariance of the problem, which might be violated more strongly in this test case compared to the circular isotropic problem, this does not come as a large surprise.
    The coarsenings depicted in~\cref{fig:corsening-aniso} demonstrate the capability of the approach to detect the anisotropy and construct coarsenings that are suitable to construct an efficient two-grid method. Even though the resulting two-grid methods for the model based approaches yield worse asymptotic convergence rates for the non grid-aligned anisotropy in c-aniso this cannot be traced to a defect in the coarsening structure. The respective sections of the coarsenings show clearly that the approach detected the anisotropy, which leads us to believe that the poor performance is more due to a poor choice of interpolatory set or interpolation weights or both rather than a poor coarsening structure.
    
    \section{Conclusion and Outlook}\label{sec:conclusion}
    In this paper we presented a new approach to adaptive algebraic multigrid construction using ideas from geostatistics and statistical learning theory. Based on the resemblance of algebraically smooth error to instances of spatial Gaussian fields we develop an empirical and semivariogram based approach to recover the covariance structure of the unknown, underlying Gaussian process. Once the covariance structure is known, efficient interpolation can be formulated by Kriging interpolation. Further exploiting the interpretation of the local interpolation error as the variance of the representation, we are able to formulate a coarsening approach that is seamlessly integrated into the determination of interpolation. Finally, by using graph distance and assuming shift invariance of the operator with respect to this pseudo-distance we are able to obtain good statistics for the semivariogram fit using only a single test vector. This is due to the fact that the semivariogram collects information at all variables to formulate a covariance function that depends solely on the distance between variables. Combined with the observation that the correlation distance of the underlying Gaussian process is very small, a single vector provides enough information about the short range correlation of values. In addition, the short correlation distance allows us to strictly localize all calculations which preserves the linear complexity of the whole process. 
    
    One apparent aim of future development is the integration of the Gaussian process approach into a multigrid setup and potentially a bootstrap type setup that is able to generate additional information about the underlying process on coarse scales as well. In line with~\cite{owhadi2017multigrid} we plan to investigate the connection between the partial differential operator, its discretization and the smoothing scheme with the resulting covariance structure of the Gaussian fields. To some extend this development can be seen in line with the investigation of optimal interpolation in algebraic multigrid methods in~\cite{BranCaoKahlFalgHu2018}, where an explicit influence of the smoother on the optimal construction of interpolation has been shown. Insight into this might allow us to translate the demonstrated potential for efficient adaptive algebraic multigrid constructions using a minimal amount of test vectors to more complex problems.
    

    \appendix
    \section{Least squares and Kriging interpolation.}\label{app:kriging-ls}
    In this appendix we present a direct comparison between the least squares~\cite{BranBranKahlLivs2011} and the Kirging interpolation for the case where $\mu=0$. As explained in Subsection \ref{sec:linear}, an alternative to estimating the stationary mean $\mu$ from the data, one can also assume $\mu=0$ if the test vectors are sampled from a centred distribution. In this case, equations \eqref{eq:average} and \eqref{eq:empCov} simplify and we obtain 
    \[
    C_{ij} = \frac{1}{K}\innerprod{V_i}{V_{j}} = \frac{1}{K}\left(V\cdot V^{T}\right)_{ij}.
    \] Correspondingly define the correlation matrix by
    \[
    X_{ij} = \frac{C_{ij}}{\sqrt{C_{ii}C_{jj}}} = \frac{\innerprod{V_i}{V_j}}{\norm{V_i}\norm{V_{j}}}.
    \]
    
    Now consider the following definition of algebraic coupling strength $\mu_{ij}$ between variables $i$ and $j$, based on the notion of least squares interpolation (cf.~\cite{BranBranKahlLivs2015b,BranBranKahlLivs2011,LivnBran2012}),
    \[
    \sigma_{ij}^2 = \frac{\sum_{k}\left(v_i^{(k)} - p_{ij} v_{j}^{(k)}\right)^2}{\sum_{k}\left(v_i^{(k)}\right)^2}.
    \] using the short-hand notation this can be written as
    \[
    \sigma_{ij}^2 = \frac{\norm{V_i - p_{ij} V_{j}}^2}{\norm{V_i}^2}.
    \] Then the minimizing $p_{ij}$ is given by
    \[
    p^{\sharp}_{ij} = \frac{\innerprod{V_i}{V_j}}{\norm{V_{j}}} = C_{ij}\cdot C_{jj}^{-1}
    \] and the corresponding minimal value by 
    \[
    (\sigma^{\sharp}_{ij})^2 = 1 - \frac{|C_{ij}|^2}{C_{ii}C_{jj}} = 1 - X_{ij}^2.
    \]
    Thus the notion of algebraic distances yields a distance with an interpretation that is similar to our approach, i.e., 
    \[
    \widehat{\sigma}_{ij}^2 = \frac{1}{1 + {K}|C_{ij}|},
    \] yields large distances for small correlations and small distances for large correlations.
    
    Based on these findings it might be worthwhile to take a closer look at the multiple interpolation variables case. Again starting with the least-squares interpolation setting and comparing it to the Kriging approach.
    
    In general the least squares interpolation for variable $i$ from variables in ${C}_i = \{j_{1},\ldots,j_{\ell}\}$ is defined as the minimizer of
    \[
    \sigma_{i,{C}_{i}}^2 = \sum_{k}\left(v_i^{(k)} - \sum_{j \in {C}_i} p_{ij} v_{j}^{(k)}\right)^2 = \sum_k \left(v_{i}^{(k)} - p_{i,{C}_{i}} V_{{C}_{i}}\right)^2.
    \]  Again using the short-hand notation with $V_{{C}_i} \in \mathbb{R}^{\ell \times K}$ this can be written as
    \[
    \sigma_{i,{C}_{i}}^2 = \norm{V_i - p_{i,{C}_i}V_{{C}_i}}^2.
    \] 
    With $C_{{C}_i,{C}_i} := V_{{C}_i} \cdot V_{{C}_i}^{T}$ we obtain the minimizer
    \[
    p_{i,{C}_i}^{\sharp} = C_{i,{C}_i} C_{{C}_i,{C}_i}^{-1}
    \]
    Correspondingly, the minimal value of $\mu_{i,{C}_i}$ is then given by
    \[
    \sigma_{i,{C}_i}^{\sharp} = C_{i,i} - C_{i,{C}_i} C_{{C}_i,{C}_i}^{-1} C_{{C}_i, i},
    \] the Schur complement of the correlation matrix.
    
    Interestingly, this last quantity is, up to appropriate scaline, equal to the estimate of the variance in our Kriging formulation. The interpolation weights of the Kriging formulation can be obtained from $p_{i,{C}_i}^{\sharp}$ by enforcing $\sum_{j \in {C}_i} p_{i,j} = 1$. That is, they are given by
    \begin{align*}
        p_{i,{C}_i}^{\rm krig} & = C_{i,{C}_i} C_{{C}_i,{C}_i}^{-1} + \frac{1 - C_{i,{C}_i} C_{{C}_i,{C}_i}^{-1}\mathbbm{1}^{T}}{\mathbbm{1}C_{{C}_i,{C}_i}^{-1}\mathbbm{1}^{T}} \mathbbm{1} C_{{C}_i,{C}_i}^{-1}\\ 
        & = p_{i,{C}_{i}}^{\sharp} + \frac{1 - C_{i,{C}_i} C_{{C}_i,{C}_i}^{-1}\mathbbm{1}^{T}}{\mathbbm{1}C_{{C}_i,{C}_i}^{-1}\mathbbm{1}^{T}} \mathbbm{1} C_{{C}_i,{C}_i}^{-1},
    \end{align*} where the latter part simply ensures $p_{i,{C}_{i}}^{\rm krig} \mathbbm{1}^{T} = 1$.
    
    Plugging the Kriging interpolation weights back into the least squares formulation yields a distance measure
    \[
    \sigma_{i,{C}_{i}}^{\rm krig} = \sigma_{i,{C}_{i}}^{\sharp} + \frac{\left(1 - C_{i{C}_{i}}C_{{C}_{i}{C}_{i}}^{-1}C_{{C}_{i}i}\right)^2}{\mathbbm{1}C_{{C}_{i}{C}_{i}}^{-1}\mathbbm{1}^{T}}.
    \]
    
    Similar constructions, explicitely preserving the constant vector, yet for reasons connected to the underlying PDE  can be found in the operator based bootstrap AMG approach~\cite{MantMcCoParkRuge2010}, where least squares interpolation is mixed with classical AMG constructions and inherits the preservation of constants from this approach. The connection between the correlation and covariance structure of test vectors has been used to define aggregation based interpolations in~\cite{LivnBran2012}, yet again a consistent definition of the whole coarsening process in terms of Gaussian processes is lacking.
   
    \bibliographystyle{siam}
\bibliography{paper_biber}
\end{document}

%% file: macros.tex
\newcommand{\innerprod}[3][2]{\langle #2,#3 \rangle_{#1}}
\newcommand{\norm}[2][2]{\| #2 \|_{#1}}